 \def\title{BERTINI AND HIS TWO FUNDAMENTAL THEOREMS}
 \def\author{S. L. Kleiman} \def\date{April 17, 1997}
\def\TheMagstep{\magstephalf}      
 \def\PaperSize{AFour}          
 \def\FirstPageOnRight{TRUE}
     \def\I#1{{\bf#1}}     \let\x=\times
\let\:=\colon \def\pd#1#2{{\partial#1\over\partial#2}}
\def\C#1{{\cal#1}}
\def\Id{I_\Delta} \def\Mx{{\bf m}_{\I x}}

\hyphenation{Fran-ces-co}

\def\TRUE{TRUE} 
\ifx\DoublepageOutput\TRUE \def\TheMagstep{\magstephalf} \fi
\mag=\TheMagstep

\parskip=0pt plus 1.75pt \parindent=2em 
\hsize=12 true cm 
\vsize=19truecm \baselineskip=\vsize \divide\baselineskip38
\abovedisplayskip6pt plus6pt minus2pt
\belowdisplayskip6pt plus6pt minus3pt

\def\today{\ifcase\month\or     
 January\or February\or March\or April\or May\or June\or
 July\or August\or September\or October\or November\or December\fi
 \space\number\day, \number\year}
\nopagenumbers
\headline={\eightpoint
  \ifnum\pageno=1\firstheadline
  \else
    \ifodd\pageno\oddheadline
    \else\evenheadline\fi
  \fi
}

\newskip\vadjustskip \vadjustskip=0.5\normalbaselineskip
\def\centertext
 {\hoffset=\pgwidth \advance\hoffset-\hsize
  \advance\hoffset-2truein \divide\hoffset by 2\relax
  \voffset=\pgheight \advance\voffset-\vsize
  \advance\voffset-2truein \divide\voffset by 2\relax
  \advance\voffset\vadjustskip
 }
\newdimen\pgwidth\newdimen\pgheight
\def\letter{letter}\def\AFour{AFour}
\ifx\PaperSize\letter
 \pgwidth=8.5truein \pgheight=11truein
 \message{- Got a paper size of letter.  }\centertext
\fi
\ifx\PaperSize\AFour
 \pgwidth=210truemm \pgheight=297truemm
 \message{- Got a paper size of AFour.  }\centertext
\fi

\def\(#1){{\rm(#1)}}

 \font\smc=cmcsc10              
\catcode`\@=11          
\def\eightpoint{\eightpointfonts
 \setbox\strutbox\hbox{\vrule height7\p@ depth2\p@ width\z@}%
 \eightpointparameters\eightpointfamilies
 \normalbaselines\rm
 }
\def\eightpointparameters{%
 \normalbaselineskip9\p@
 \abovedisplayskip9\p@ plus2.4\p@ minus6.2\p@
 \belowdisplayskip9\p@ plus2.4\p@ minus6.2\p@
 \abovedisplayshortskip\z@ plus2.4\p@
 \belowdisplayshortskip5.6\p@ plus2.4\p@ minus3.2\p@
 }
\newfam\smcfam
\def\eightpointfonts{%
 \font\eightrm=cmr8 \font\sixrm=cmr6
 \font\eightbf=cmbx8 \font\sixbf=cmbx6
 \font\eightit=cmti8
 \font\eightsmc=cmcsc8
 \font\eighti=cmmi8 \font\sixi=cmmi6
 \font\eightsy=cmsy8 \font\sixsy=cmsy6
 \font\eightsl=cmsl8 \font\eighttt=cmtt8}
\def\eightpointfamilies{%
 \textfont\z@\eightrm \scriptfont\z@\sixrm  \scriptscriptfont\z@\fiverm
 \textfont\@ne\eighti \scriptfont\@ne\sixi  \scriptscriptfont\@ne\fivei
 \textfont\tw@\eightsy \scriptfont\tw@\sixsy \scriptscriptfont\tw@\fivesy
 \textfont\thr@@\tenex \scriptfont\thr@@\tenex\scriptscriptfont\thr@@\tenex
 \textfont\itfam\eightit        \def\it{\fam\itfam\eightit}%
 \textfont\slfam\eightsl        \def\sl{\fam\slfam\eightsl}%
 \textfont\ttfam\eighttt        \def\tt{\fam\ttfam\eighttt}%
 \textfont\smcfam\eightsmc      \def\smc{\fam\smcfam\eightsmc}%
 \textfont\bffam\eightbf \scriptfont\bffam\sixbf
   \scriptscriptfont\bffam\fivebf       \def\bf{\fam\bffam\eightbf}%
 \def\rm{\fam0\eightrm}%
 }
\def\vfootnote#1{\insert\footins\bgroup
 \eightpoint 
 \interlinepenalty\interfootnotelinepenalty
  \splittopskip\ht\strutbox 
  \splitmaxdepth\dp\strutbox \floatingpenalty\@MM
  \leftskip\z@skip \rightskip\z@skip \spaceskip\z@skip \xspaceskip\z@skip
  \textindent{#1}\footstrut\futurelet\next\fo@t}

 \newdimen\fullhsize \newbox\leftcolumn
 \def\fulline{\hbox to \fullhsize}
\def\doublepageoutput
{\let\lr=L
 \output={\if L\lr
          \global\setbox\leftcolumn=\columnbox \global\let\lr=R%
        \else \doubleformat \global\let\lr=L\fi
        \ifnum\outputpenalty>-20000 \else\dosupereject\fi}%
 \def\doubleformat{\shipout\vbox{%
        \fulline{\hfil\hfil\box\leftcolumn\hfil\columnbox\hfil\hfil}%
                                }%
                  }%
 \def\columnbox{\vbox
   {\makeheadline\pagebody\makefootline\advancepageno}%
   }%
 \fullhsize=\pgheight \hoffset=-1truein
 \voffset=\pgwidth \advance\voffset-\vsize
  \advance\voffset-2truein \divide\voffset by 2
  \advance\voffset\vadjustskip
 \let\firstheadline=\hfil
 
\ifx\FirstPageOnRight\TRUE 
 \null\vfill\nopagenumbers\eject\pageno=1\relax
\fi
}
\ifx\DoublepageOutput\TRUE \doublepageoutput \fi

\def\sctvskip{3\normalbaselineskip}
\def\sctvspace{\removelastskip\vskip0pt plus\sctvskip
 \penalty-250 \vskip0pt plus-\sctvskip \bigskip\medskip
}
\def\sct#1\par
 {\sctvspace \advance\sctno by1\relax
 {\bf \number\sctno.~#1\unskip.}
 \medskip
}
\newcount\sctno \sctno=0
\def\sctnum#1{\sctno=#1\advance\sctno by -1}

\def\art#1#2{\medbreak{\bf#1 (\number\sctno.#2)}\enspace}
\def\proclaim#1 #2{\art{#1}{#2}\begingroup\it}

\def\endproclaim{\endgroup\medbreak}
\def\Cs#1){(\number\sctno.#1)}
\def\part#1 {\par{\rm (#1)\enspace}\ignorespaces}

\def\fig#1{Fig.\penalty\@M \thinspace#1}

 \newcount\refno \refno=0        \def\NoKey{*!*}
 \def\MakeKey{\advance\refno by 1 \expandafter\xdef
  \csname\TheKey\endcsname{{\number\refno}}\NextKey}
 \def\NextKey#1 {\def\TheKey{#1}\ifx\TheKey\NoKey\let\next\relax
  \else\let\next\MakeKey \fi \next}
 \def\RefKeys #1\endRefKeys{\expandafter\NextKey #1 *!* }
 \def\SetRef#1 #2,{\hang
        \llap{[\csname#1\endcsname]\enspace}%
        #2,%
}
 \newbox\keybox \setbox\keybox=\hbox{[8]\enspace}
 \newdimen\keyindent \keyindent=\wd\keybox
\def\references{
  \bgroup   \frenchspacing   \eightpoint
   \parindent=\keyindent  \parskip=\smallskipamount
   \everypar={\SetRef}\par
}
\def\endreferences{\egroup}
 \def\ms#1 {Mem. Scelte, Zanichelli, Bologna, 1956--66, \vol.#1\unskip, }
 \def\cp#1 {Collected Papers, MIT Press, 1972--73, {\bf #1}, }
 \def\serial#1#2{\expandafter\def\csname#1\endcsname ##1 ##2
        ##3\par{\unskip\ #2 {\bf##1} (##2), ##3\unskip.\par}} 
 \serial{ajm}{Amer. J. Math.}
 \serial{adm}{Ann. Mat. Pura Appl.}
 \serial{am}{Ann. of Math.}
 \serial{ahes}{Arch. History Exact Sci.}
 \serial{aast}{Atti Accad. Sci. Torino Cl. Sci. Fis. Mat. Natur.}
 \serial{bdmf}{Boll. di mat. Firenze}
 \serial{bumi}{Boll. Un. Mat. Ital.}
 \serial{emcmc}{Esercitazioni matematiche (Circ. mat. Catania)}
 \serial{jram}{J. Reine Angew Math.}
 \serial{ma}{Math. Ann.}
 \serial{mast}{Mem. Accad. Sci. Torino Cl. Sci. Fis. Mat. Natur.}
 \serial{mcsuk}{Mem. Coll. Sci. Univ. Kyoto Ser. A. Math.}
 \serial{rcmp}{Rend. Circ. Mat. Palermo}
 \serial{xlm}{Rend. Accad. Naz. Sci. XL Mem. Mat.}
 \serial{anlr}{Atti Accad. Naz. Lincei Cl. Sci. Fis. Mat. Natur. Rend. Lincei}
 \serial{rral}{Atti Accad. Naz. Lincei Cl. Sci. Fis. Mat. Natur. Mem.}
 \serial{rril}{Istit. Lombardo Accad. Sci. Lett. Rend. A Istituto}
 \serial{tams}{Trans. Amer. Math. Soc.}

\def\UThin{\penalty\@M \thinspace\ignorespaces}
\def\p.{p.\UThin} \def\pp.{pp.\UThin} \def\vol.{vol.\UThin}
\def\relaxnext@{\let\next\relax}
\def\cite#1{\relaxnext@
 \def\nextiii@##1,##2\end@{\unskip\space
  {\let~=\UThin\rm[\SetKey{##1},##2]}}%
 \in@,{#1}\ifin@\def\next{\nextiii@#1\end@}\else
 \def\next{{\rm[\SetKey{#1}]}}\fi\next}
\newif\ifin@
\def\in@#1#2{\def\in@@##1#1##2##3\in@@
 {\ifx\in@##2\in@false\else\in@true\fi}%
 \in@@#2#1\in@\in@@}
\def\SetKey#1{{\bf\csname#1\endcsname}}

\catcode`\@=12 

\chardef\atcode=\catcode`\@
\catcode`\@=11

\def\PSpicture#1#2{\vbox to #2\bgroup\gdef\PSstr@ng{"}%
  \offinterlineskip
  \hrule height0pt width #1\relax\vfill\PSpictsetup}
\def\endPSpicture{\special{\PSstr@ng}\egroup}

\newdimen\PSunit \PSunit=1bp

\def\llcoords#1{\begingroup\c@l=0 \st@ck#1 \end
  \ifnum\c@l=2 \else
      \errmessage{Two coordinates needed for \string\llcoords}\fi
   \dimen@x=\ll@x\dimen@y=\ll@y
   \global\ll@y=-\Acol \advance\c@l-1 \global\ll@x=-\Acol
   \advance\ll@x-\dimen@x \advance\ll@y-\dimen@y
   \PS{\The\ll@x\space\The\ll@y\space translate}\endgroup}

\newdimen\PSlinewidth \PSlinewidth=0.5pt

\def\setPSlinewidth#1{\PSlinewidth=#1 %
   \PS{\The\PSlinewidth\space\The\PSunit\space div setlinewidth}}

\def\PS#1{\global\let\tmp=\PSstr@ng
   \global\edef\PSstr@ng{\tmp\space#1}}

\def\PSpictsetup{\newlocregs \ll@x=0pt \ll@y=0pt %
  \count@=\PSunit\PS{\the\count@\space 65781.76 div dup scale}%
  \setPSlinewidth\PSlinewidth
  \defPSarrowh
  \setPSarrowheadlength{10\PSlinewidth}}

\def\TX#1#2{\begingroup\c@l=0 \st@ck#2 \end
  \ifnum\c@l=2 \else
     \errmessage{Two coordinates needed for \string\TX}\fi
  \advance\ll@y\Acol \advance\c@l-1 \advance\ll@x\Acol
  \ll@y=\The\ll@y\PSunit \ll@x=\The\ll@x\PSunit
  \vbox to 0pt{\kern-\ll@y\vss
       \hbox to 0pt{\kern\ll@x#1\kern-\ll@x}
       \vss\kern\ll@y}\endgroup}

\def\markbox{\hss\vrule height 3\PSlinewidth
      width 3\PSlinewidth \hss}

\newcount\base \base=32768

\def\m@ltf#1#2{\dimen@y=#1\divide\dimen@y\base
  \dimen@x=-\dimen@y \multiply\dimen@x\base\advance\dimen@x#1%
  \multiply\dimen@y#2\multiply\dimen@x#2\divide\dimen@x\base
  \advance\dimen@y\dimen@x#1=\dimen@y}

\def\inv@rt#1#2{\count@@=#1\multiply\count@@\base
  \advance\count@@-#2#2=\base \multiply#2\base \divide#2\count@@}

\def\newlocdimen#1{\advance\count11 by \@ne
  \ifnum\count11 <\insc@unt
     \else\errmessage{No room for the spline}\fi
  \dimendef#1=\count11 }
\def\newloccount#1{\advance\count10 by \@ne
  \ifnum\count10 <\insc@unt
     \else\errmessage{No room for the splines}\fi
  \countdef#1=\count10 }
\def\newdimenvar#1{\newlocdimen\currvar
   \expandafter\let\csname#1\the\c@l\endcsname=\currvar
   \wlog{\expandafter\string\csname#1\the\c@l\endcsname
      =\string\dimen\the\count11 }}
\def\newcountvar#1{\newloccount\currvar
   \expandafter\let\csname#1\the\c@l\endcsname=\currvar
  \wlog{\expandafter\string\csname#1\the\c@l\endcsname
   =\string\count\the\count10 }}

\def\newlocregs{\newloccount\c@l\newloccount\numc@ls
  \newloccount\prevF\newloccount\prevQ
  \newloccount\count@\newloccount\count@@
  \newlocdimen\dimen@\newlocdimen\dimen@x\newlocdimen\dimen@y
  \newlocdimen\ll@x\newlocdimen\ll@y}

\def\Acol{\csname A\the\c@l\endcsname}
\def\Fcol{\csname F\the\c@l\endcsname}
\def\Qcol{\csname Q\the\c@l\endcsname}
\def\Zcol{\csname Z\the\c@l\endcsname}

\def\spline{\begingroup\spl@ne}
\def\cspline{\begingroup\let\d@spline=\d@cspline
   \let\PSw@finish=\PSw@cfinish\spl@ne}
\def\spl@ne#1{\maybem@rks\c@l=0 \st@ck#1 \end
   \ch@ckeven\ch@ckthree
   \d@spline\d@PSw\PSw@finish\stroke\endgroup}
\def\ch@ckthree{\ifnum\numc@ls<6 %
      \errmessage{At lest three points needed for SPLINE}\fi}
\def\ch@cktwo{\ifnum\numc@ls<4 %
      \errmessage{At lest two points needed}\fi}
\def\ch@ckeven{\ifodd\numc@ls
       \errmessage{Even number of coordinates needed}\fi}
\def\d@spline{\fwd@init\loop\fwd@body\ifnum\c@l<\numc@ls\repeat
  \fwd@finish\s@lve}
\def\d@cspline{\fwd@cinit\loop\fwd@cbody\ifnum\c@l<\numc@ls\repeat
  \fwd@cfinish\cs@lve}

{\def\\{\global\let\sptoken= }\\ } 
\def\st@ck{\futurelet\nexttok\testst@ck}
\def\testst@ck{\ifx\nexttok\sptoken
   \let\next\eatspacest@ck
   \else\ifx\nexttok\end\let\next\stackfinish
   \else\let\next\st@ckit\fi\fi\next}
\def\eatspacest@ck#1 {\st@ck}
\def\stackfinish#1{\numc@ls=\c@l}
\def\st@ckit#1 {\advance\c@l1 \newdimenvar A\Acol=#1pt %
  \markoptions\st@ck}

\def\fwd@init{
  \c@l=1 \newdimenvar Z\currvar=\Acol
    \advance\c@l2 \advance\currvar 2\Acol
  \c@l=2 \prevF=\base \divide\prevF2 \newcountvar F\Fcol=\prevF
    \newdimenvar Z\currvar=\Acol
    \advance\c@l2 \advance\currvar 2\Acol
  \c@l=4 }
\def\fwd@cinit{%
  \c@l=\numc@ls
   \newcountvar F\let\lastF=\currvar \lastF\base \multiply\lastF4 %
   \newdimenvar Z\let\lastY=\currvar \lastY=4\Acol
  \advance\c@l-1 \newdimenvar Z\let\lastX=\currvar\lastX=4\Acol
  \c@l=1 \advance\lastX 2\Acol
   \prevQ=\base\newcountvar Q\Qcol=\prevQ
   \newdimenvar Z\Zcol=4\Acol\advance\c@l2 \advance\currvar 2\Acol
  \c@l=2 \advance\lastY 2\Acol
   \prevF=\base \divide\prevF4 \newcountvar F\Fcol=\prevF
   \newdimenvar Z\Zcol=4\Acol\advance\c@l2 \advance\currvar 2\Acol
  \c@l=4 }

\def\fwd@body{%
  \advance\c@l-1 \newdimenvar Z%
  \advance\c@l-2 \currvar=-\Zcol \m@ltf\currvar\prevF
  \advance\c@l2 \advance\currvar 4\Acol
  \advance\c@l2 \advance\currvar 2\Acol
  \advance\c@l-2 \ifodd\c@l\else
     \inv@rt4\prevF\newcountvar F\Fcol=\prevF\fi
  \advance\c@l2 }
\def\fwd@cbody{%
  \advance\c@l-1 \newdimenvar Z%
  \advance\c@l-2\currvar=-\Zcol \m@ltf\currvar\prevF
  \advance\c@l2 \advance\currvar 4\Acol
  \advance\c@l2 \advance\currvar 2\Acol
  \advance\c@l-2 \ifodd\c@l\else
     \count@=\prevQ
     \prevQ=-\prevQ \multiply\prevQ\prevF \divide\prevQ\base
     \multiply\count@\prevQ \divide\count@\base
     \advance\lastF\count@
     \advance\c@l-2 \dimen@=\Zcol \m@ltf\dimen@\prevQ
       \advance\lastY\dimen@
     \advance\c@l-1 \dimen@=\Zcol \m@ltf\dimen@\prevQ
       \advance\lastX\dimen@
     \advance\c@l2 \newcountvar Q\Qcol=\prevQ
     \advance\c@l1 \inv@rt4\prevF\newcountvar F\Fcol=\prevF\fi
  \advance\c@l2 }

\def\fwd@finish{%
  \c@l=\numc@ls
  \advance\c@l-1 \newdimenvar Z%
  \advance\c@l-2 \currvar=-\Zcol\m@ltf\currvar\prevF
  \advance\c@l2 \advance\currvar 3\Acol
  \advance\c@l1 \newdimenvar Z%
  \advance\c@l-2 \currvar=-\Zcol\m@ltf\currvar\prevF
  \advance\c@l2 \advance\currvar 3\Acol
  \inv@rt2\prevF \newcountvar F\Fcol=\prevF}
\def\fwd@cfinish{%
  \c@l=\numc@ls \advance\prevQ\base
    \count@=-\prevQ \multiply\count@\prevF\divide\count@\base
  \advance\c@l-2 \dimen@=\Zcol
    \m@ltf\dimen@\count@\advance\lastY\dimen@
  \advance\c@l-1 \dimen@=\Zcol
    \m@ltf\dimen@\count@\advance\lastX\dimen@
    \multiply\count@\prevQ \divide\count@\base\advance\lastF\count@
  \prevF=\base \multiply\prevF\base \divide\prevF\lastF \lastF=\prevF
  \c@l=\numc@ls}

\def\s@lve{\c@l=\numc@ls \m@ltf\Zcol\prevF
  \advance\c@l-1 \m@ltf\Zcol\prevF
  \loop\advance\c@l-1 %
    \ifodd\c@l \else \prevF=\Fcol\fi
    \advance\c@l2 \dimen@=-\Zcol
    \advance\c@l-2 \advance\Zcol\dimen@\m@ltf\Zcol\prevF
    \ifnum\c@l>1 \repeat}
\def\cs@lve{\c@l=\numc@ls \m@ltf\Zcol\prevF
  \advance\c@l-1 \m@ltf\Zcol\prevF
  \loop\advance\c@l-1 %
    \ifodd\c@l \dimen@=-\lastX
         \else \dimen@=-\lastY \prevF=\Fcol
               \advance\c@l-1 \prevQ=\Qcol\advance\c@l1 \fi
    \m@ltf\dimen@\prevQ
    \advance\c@l2 \advance\dimen@-\Zcol
    \advance\c@l-2 \advance\dimen@\Zcol \m@ltf\dimen@\prevF
    \Zcol=\dimen@ \ifnum\c@l>1\repeat}

{\catcode`p=12 \catcode`t=12 \gdef\STRIPPT#1pt{#1}}
\def\The#1{\expandafter\STRIPPT\the#1}

\def\PSw@#1{\advance\c@l-1 \PS{\The#1}\advance\c@l1 \PS{\The#1}}
\def\PSw@Wcol{\advance\c@l1 \dimen@x=2\Acol\advance\dimen@x-\Zcol
  \advance\c@l1 \dimen@y=2\Acol\advance\dimen@y-\Zcol
  \PS{\The\dimen@x\space\The\dimen@y}}
\def\PSw@init{\c@l=2 \newpath\PSw@\Acol\moveto}
\def\PSw@ZWA{%
  \PSw@\Zcol\PSw@Wcol\PSw@\Acol\PS{curveto}}
\def\PSw@finish{\relax}
\def\PSw@cfinish{\PSw@\Zcol
  \c@l=0 \PSw@Wcol\PSw@\Acol\PS{curveto}}
\def\d@PSw{\PSw@init
  \loop\PSw@ZWA\ifnum\c@l<\numc@ls\repeat}
\def\normalpaths{\def\stroke{\PS{stroke}}%
  \def\moveto{\PS{moveto}}\def\newpath{\PS{newpath}}}
\normalpaths

\def\PSline{\begingroup\PSl@ne}
\def\PSl@ne#1{\maybem@rks\c@l=0 \st@ck#1 \end
  \ch@ckeven\ch@cktwo
  \c@l=2 \newpath\PSw@\Acol\moveto
  \loop\advance\c@l2 \PSw@\Acol\PS{lineto}%
  \ifnum\c@l<\numc@ls\repeat\PSl@nefinish\stroke\endgroup}
\let\PSl@nefinish=\relax
\def\PScline{\begingroup\def\PSl@nefinish{\PS{closepath}}\PSl@ne}

\def\PSarc{\begingroup\PS@rc}
\def\PSarcn{\begingroup\let\PS@arcfinish=\PS@arcnfinish\PS@rc}
\def\PS@rc#1#2{\c@l=0 \st@ck#2 #1 \end
  \ifnum\c@l=5 \else
     \errmessage{Wrong number of arguments for circles/arcs}\fi
  \newpath\c@l=0 \loop
       \advance\c@l1 \PS{\The\Acol}\ifnum\c@l<5 \repeat
  \PS@arcfinish\stroke\endgroup}
\def\PS@arcfinish{\PS{arc}}
\def\PS@arcnfinish{\PS{arcn}}

\let\markoptions=\relax
\def\putTeXmarks{\def\maybem@rks{\let\markoptions=\TeXm@rks}}
\def\putPSmarks{\def\maybem@rks{\let\markoptions=\PSm@rks}}
\def\nomarks{\let\maybem@rks=\relax}

\nomarks
\def\TeXm@rks{\ifodd\c@l\else
  \begingroup \advance\ll@y\Acol \advance\c@l-1 \advance\ll@x\Acol
  \ll@y=\The\ll@y\PSunit \ll@x=\The\ll@x\PSunit
  \vbox to 0pt{\kern-\ll@y\vss
       \hbox to 0pt{\kern\ll@x\markbox\kern-\ll@x}
       \vss\kern\ll@y}\endgroup\fi}
\let\PSmarkradius=\PSlinewidth
\def\PSm@rks{\ifodd\c@l\else\PS{newpath}\PSw@\Acol
   \PS{\The\PSmarkradius\space\The\PSunit\space div
         0 360 arc stroke}\fi}
\def\TeXmark{\begingroup\let\markoptions=\TeXm@rks\m@rk}
\def\PSmark{\begingroup\let\markoptions=\PSm@rks\m@rk}
\def\m@rk#1{\c@l=0 \st@ck#1 \end\ch@ckeven\endgroup}

\def\setdash#1{\dimen@=#1\relax
  \PS{[\The\dimen@\space\The\PSunit\space div] 0 setdash}}
\def\nodashes{\PS{[] 0 setdash}}

\def\longpath#1{\c@l=0 \st@ck#1 \end
  \ifnum\c@l=2 \else
     \errmessage{Wrong number of arguments for \string\longpath}\fi
  \PSw@\Acol\PS{moveto}%
  \let\newpath\relax\let\stroke\relax\def\moveto{\PS{lineto}}}

\def\subpaths{\PS{newpath}\let\newpath\relax
   \let\stroke\relax\def\moveto{\PS{moveto}}}

\def\defPSarrowh{\PS{/arrowh {
 /curmtx matrix currentmatrix def
 currentpoint translate
 arrlen dup scale
 rotate
 -1  0
 -2.5  0.5
 -3  1
 curveto
 -2.5  0.5
 -2.5  -0.5
 -3  -1
 curveto
 -2.5  -0.5
 -1  0
 0 0
 curveto closepath fill
 curmtx setmatrix} def}}
\def\setPSarrowheadlength#1{\dimen@=#1\relax
  \PS{/arrlen \The\dimen@\space\The\PSunit\space div 3 div def}}

\def\PSarrow#1{\begingroup
  \c@l=0 \st@ck#1 \end
  \ifnum\c@l=4 \else
  \errmessage{4 coordinates needed for \string\PSarrow}\fi
  \PS{newpath}\c@l=2 \PSw@\Acol\PS{moveto}\c@l=4 \PSw@\Acol
  \PS{lineto}\stroke\PS{newpath}\PSw@\Acol\PS{moveto}
  \dimen@y=\Acol\advance\c@l-1 \dimen@x=\Acol
  \advance\c@l-1 \advance\dimen@y-\Acol
  \advance\c@l-1 \advance\dimen@x-\Acol
  \ifdim\dimen@x=0pt\ifdim\dimen@y=0pt
      \errmessage{equal endpoints of \string\PSarrow}\fi\fi
  \PS{\The\dimen@y\space\The\dimen@x\space atan arrowh}\endgroup}

\catcode`\@=\atcode
\RefKeys
 Bert77 Bert82 Bert01 Bert04 BertI BertC Berz15 Berz33a Berz33b Bot77
Bot81 BoCoGa96 BrCi95 BrCiSe92 Bri92 Cast33 Con34 Cool Dieu En93 En96
En38 EnCa32 EnCh18 Fub33 Greit Lu93 Lu94 Mat51 Mum72 Noe73 Noe75 Par91
PiSi06 Poin Ros71 Scor34 Sev06 Sev08 vdW35 vdW37 Ver03 Volt Zar35 Zar41
Zar44 Zar44r Zar47 Zar58
  \endRefKeys

 \let\headtitle=\title
\def\firstheadline{\vtop{\eightrm
 \hbox{RENDICONTI DEL CIRCOLO MATEMATICO DI PALERMO}
 \hbox{Serie II - Supplemento (1997) $=$ alg-geom/9704018}
 }\hfill}
 \let\headtitle=\title
\def\oddheadline{\hfil\headtitle\hfil\llap{\folio}}
\def\evenheadline{\rlap{\folio}
 \hfil\author\hfil\llap{\date}}

\leavevmode
 {\bigskipamount=8\bigskipamount \bigskip}
 \centerline{\bf\title}
 \bigskip
 \centerline{STEVEN L. KLEIMAN}
 \bigskip\bigskip
{\narrower\eightpoint\noindent{\bf Abstract.}
 After reviewing Bertini's life story, a fascinating drama, we make a
critical examination of the old statements and proofs of Bertini's two
fundamental theorems, the theorem on variable singular points and the
theorem on reducible linear systems.  We explain the content of the
statements in a way that is accessible to a nonspecialist, and we
develop versions of the old proofs that are complete and rigorous by
current standards.  In particular, we prove a new extension of Bertini's
first theorem, which treats variable $r$-fold points for any $r$.
 \par}
 \bigskip

\sct Preface

Eugenio Bertini (1846--1933) studied in the 1860s with Luigi Cremona
(1830--1903), the father of Italian algebraic geometry.  Bertini was
Cremona's first student, and one of his best.  At the time,
Cremona was developing the first general theory of birational
transformations of the plane and of 3-space, the transformations that
now bear his name.  Bertini was attracted to the subject and advanced
it; thus he was led to prove, in his paper \cite{Bert82} dated December
1880 and published in 1882, the two theorems that now bear his name:
Bertini's theorem on variable singular points, and Bertini's theorem
on reducible linear systems.

Bertini's two theorems soon became fundamental tools in algebraic
geometry.  In his comprehensive treatise of 1931 on plane curves
\cite{Cool}, Julian Coolidge (1873--1954) wrote on \p.115 that Bertini's
first theorem ``will be of utmost importance to us,'' and then called it
``vitally important'' on the next page.  In his obituary \cite{Cast33,
p.~747} of Bertini, Guido Castelnuovo (1865--1952) wrote that both
theorems ``come into play at every step in all the work'' of the Italian
school of algebraic geometry.  In his obituaries, \cite{Berz33a, p.~621}
and \cite{Berz33b, p.~150}, Luigi Berzolari (1863--1949) wrote that the
theorems are ``now classical and of constant use in current research in
algebraic geometry.''

Bertini's theorems are today as fundamental as ever in algebraic
geometry.  Moreover, from about 1880 until 1950, their statements and
proofs evolved in generality, rigor, and clarity in the hands notably of
Bertini himself, of Frederigo Enriques (1871--1946), of Bartel van der
Waerden (1903--96), and of Oscar Zariski (1899--1986).  For these reasons
alone, the theorems make an appropriate subject for a historical
analysis.  However, in addition, it turns out that the old statements
and proofs are rather interesting from a purely mathematical point of
view.  For example, Bertini's original first theorem treated the
variable $r$-fold points for an arbitrary $r$, not simply the variable
singular points.  Furthermore, Bertini's second theorem used to be
derived from the first.  Yet the old statements and proofs have been
nearly forgotten, replaced by various new ones.

The bulk of the present article is devoted to a critical examination of
the old statements and proofs.  Section~3 discusses the first theorem
for an ambient projective space over the complex numbers, the case that
Bertini considered originally.  Section~4 discusses the extension of
the first theorem to an arbitrary ambient variety; in particular, it
proves a new general version, Theorem~(4.4), for an arbitrary $r$.
Finally, Section~5 discusses the second theorem.  Nothing is said
anywhere about the newer theory, despite its importance in contemporary
algebraic geometry.  On the other hand, an attempt is made throughout
to develop versions of the old statements and proofs that are
acceptably complete, rigorous, and clear by current standards.  At the
same time, an attempt is made to explain the content of the statements
and the spirit of the proofs in a way that is accessible to a
nonspecialist.

The life stories of mathematicians are often fascinating human dramas,
which show the enormous influence that social circumstances have on the
development of mathematics.  Thus it is with Bertini's life story, which
we'll review in Section 2, drawing on a number of sources, including the
historical articles of Bottazzini, \cite{Bot77} and \cite{Bot81}, and of
Vito Volterra (1860--1940) \cite{Volt}, the historical monographs of
Aldo Brigaglia and Ciro Ciliberto \cite{BrCi95} and of Enriques
\cite{En38, \pp.281--92}, and the announcement \cite{Bri92} of the death
of Enrico Betti (1823--92) made by Francesco Brioschi (1824--97).  We'll
also draw on the obituary of Bertini written by his close colleague
Castelnuovo \cite{Cast33} (who received 185 conserved personal letters
from him, see \cite{BoCoGa96, \p.XXIX}), the obituaries by Bertini's
student at Pavia, Berzolari, \cite{Berz33a} and \cite{Berz33b}, and the
obituaries by his students at Pisa, Alberto Conti (1873--1940)
\cite{Con34}, Guido Fubini (1879--1943) \cite{Fub33} and Gaetano Scorza
(1876--1939) \cite{Scor34}.

The obituaries carry all the authority and compassion of their authors'
firsthand knowledge of Bertini's character, his lecturing and his
writing.  Further evidence of his marvelous character is found in nearly
60 of the 668 letters from Enriques to Castelnuovo, published in
\cite{BoCoGa96}.  Each source contains its own gems, and all can be
recommended.  The obituaries also contain surveys of Bertini's
scientific work, and some of this material will be repeated briefly
below; see also the technical monographs of Berzolari \cite{Berz15,
p.~328}, Coolidge \cite{Cool, p.~481}, and Jean Dieudonne (1906--92)
\cite{Dieu, p.~111}.  Two of the obituaries, Berzolari's \cite{Berz33b}
and Fubini's \cite{Fub33}, contain Bertini's scientific bibliography; of
these two, the first is more carefully prepared, and it also contains
the fullest discussion (19 pages) of the works themselves.  However, of
Bertini's original discoveries, aside from the two fundamental theorems,
only the classification of plane involutions remains significant, and it
is of less general interest.

Because Bertini's life was so intertwined with Cremona's, the latter
will also be reviewed, although more briefly, following Greitzer's short
biography \cite{Greit} and following the obituaries written by his two
students, Giuseppe Veronese (1857--1917) \cite{Ver03} and Bertini
\cite{Bert04}.

The present article is an expanded version of a talk given at the
conference, {\it Algebra e Geometria (1860--1940): Il contributo
Italiano,} which took place in Cortona, Italy, 4--8 May 1992 under the
scientific direction of Brigaglia, Ciliberto and Edoardo Sernesi; see
the proceedings \cite{BrCiSe92}.  The author's conversations with Rick
Miranda and with Sernesi in Cortona and recently with Anders Thorup in
Copenhagen led to the proof of the new extension, Theorem~(4.4), of
Bertini's first theorem.  Beverly Kleiman provided valuable help with
the copy-editing of this article.  Its appearance, at last, is due to
the long-standing persistent encouragement of Umberto Bottazzini, who
also kindly provided copies of a fair number of hard-to-find articles.
It is a pleasure for the author to express his heartfelt thanks to all
of these individuals, especially to Umberto.

\sct Bertini's Life

Eugenio Bertini was born on 8 November 1846 in Forli, Italy, about
halfway between Rimini and Bologna, to Vicenzo Bertini, a typographer,
and his wife Agata, n\'ee Bezzi.  Bertini attended secondary school for
two years at the technical institute in Forli, where he showed distinct
aptitude in mathematics.  The family was poor, but with the financial
assistance of the Congregazione di Carit\`a di Forli, Bertini was able
to go on for a higher education.  In 1863, nearly 17, he entered the
University of Bologna, intending to become an engineer, but he was
drawn into pure mathematics by the lectures of Luigi Cremona, who was
nearly twice his age at the time.

This was the period of the unification of Italy.  There were three wars
of independence, of 1848, 1859 and 1866.  In 1859 Lombardy was annexed
by Sardinia--Piedmont, after Cavour had gotten France to help drive out
Austria.  In 1860 the northern states of Tuscany, Modena, Parma, and
Romagna rose up against their princes, and were annexed with France's
permission in return for Savoy and Nice.  Scientists, including the
mathematicians Betti and Brioschi, played prominent roles in the
enlightened new government, which promptly established chairs of higher
mathematics.  Professors and students could move more freely and easily
from university to university, as the universities were now under a
single ministry.

Just before, in the spring of 1858, Betti, Brioschi and several others
met in Genoa to found the journal {\it Annali di matematica pura e
applicata}, modeling it after {\it Crelle's Journal f\"ur die reine und
angewande Mathematik\/} and {\it Liouville's Journal de math\'ematiques
pures et applique\'ees}; a decade later under the editorship of
Brioschi and Cremona, it grew into one of the great European journals.
In the fall of 1858, Betti, Brioschi and Felice Casorati (1835--90)
visited the mathematical centers of G\"ottingen, Berlin, and Paris, and
made many important mathematical contacts.

Many Italian mathematicians did pure research of international caliber,
yet also made an effort to write good texts and to train skilled
engineers, who were needed to build the industry and the infrastructure
of the new nation.  Italy felt a particular kinship with Germany, which
was also, at the time, forging a nation out of a maze of states.  So
Italian mathematicians studied the works of Gauss, Jacobi and Riemann,
and of Pl\"ucker, Clebsch, and Noether.  In 1866 Italy joined Germany
(Prussia) in a war against Austria, which gave the Venetia to France,
which in turn gave it to Italy.  Bertini volunteered as an infantryman
under Garabaldi in this brief third war for Italian independence.

Eighteen years earlier, in April 1848, Cremona had interrupted his
studies to fight as a volunteer in the first war against Austrian rule,
attaining the rank of sergeant.  He took part in the heroic defense of
Venice, which capitulated on 24 August 1849.  Because of the discipline
and bravery of the defenders, they were allowed by the victors to leave
as a unit to serve as a model of military and civil virtue.  Cremona
returned to his native Pavia to find that his mother had died a few
months earlier.  Soon he became gravely ill with typhus, which he had
picked up in the war.  Strong willed, he reentered the university the
same year.  From his teacher Brioschi, he learned to love science and
to pass on this love in his own lectures and texts, which were
excellent.  On 9 May 1853, he was awarded his laurea degree in civil
engineering and architecture ({\it sic\/}!).

Cremona could not enter the official educational system right away
because of his record of military service against Austria.  So he became
a private tutor in Pavia for several years (a common profession at the
time).  On 22 November 1855 he was appointed provisionally as a teacher
at Pavia's ginnasio liceale, and on 17 December 1856 he was promoted to
the rank of associate.  A month later he was appointed to full rank at
the ginnasio in the city of Cremona.  He remained there for three years,
until the new government of Lombardy appointed him to the Liceo S.\
Alessandro in Milan and then on 10 June 1860 to the first chair of
higher geometry in Bologna.  There he carried out his most important
original research, which concerned birational transformations and their
applications.  For part of this work, he shared the 1864 Steiner prize
in geometry awarded by the Berlin Academy (the other winner was R.
Sturm).  In 1863, Brioschi founded an engineering school in Milan, and
on his recommendation, Cremona was transferred there in October 1866.

After the third war, Bertini was advised by Cremona to resume his
studies under Betti and Ulisse Dini (1845--1918) in Pisa.  Bertini
earned his laurea in 1867 at the Ateneo pisano and his teaching
certificate in 1868 at the Scuola Normale Superiore.  Nominated
immediately afterwards to the chair of mathematics at the Liceo Parini
in Milan, Bertini had the opportunity, the next academic year
(1868--69), to attend a three-part course given by Brioschi, Casorati,
and Cremona, on Abelian integrals from the three different points of
view: Jacobi's analysis, Riemann's topology, and Clebsch's algebraic
geometry.  This course had a profound and lasting affect on Bertini.
In particular, he was led to give, in 1869, one of the first and
simplest geometric proofs of the invariance of the genus of a curve
under birational correspondence.  This proof was included in three
standard texts (Clebsch--Lindemann of 1876, Salmon--Chemin of 1884, and
Enriques--Chisini of 1918 \cite{EnCh18, \vol.2, \pp.131-35}), making
Bertini's name famous.

In the fall of 1871, Bertini took the chair of mathematics at the Liceo
Visconti in Rome, but in addition he taught descriptive geometry as an
adjunct at the university.  In the fall of 1873, he abandoned
secondary-school teaching, and began teaching projective geometry too at
the university.  He remained in Rome one more year only.  In 1875, he
won, with Betti's support, a competition for the chair in advanced
geometry in Pisa.  He served first as an adjunct, but after three years,
was promoted to a regular position, and he remained in Pisa two more
years, until 1880.

On 9 October 1873, Cremona was appointed director of the newly
established engineering school in Rome, and he soon became so burdened
with administration that he had little time left over for creative
research.  When Bertini went to Pisa in 1875, Cremona objected
vociferously, and only forgave him two years later when Bertini,
generous as always, offered to give up his chair in favor of his
teacher, who appeared for a moment to want to leave Rome.  However, in
November 1877, Cremona was appointed to a chair at the University of
Rome.  On 16 March 1879 he was appointed a senator, and then his
research activities stopped completely.  On 10 June 1903, he left his
sickbed to act on some legislation, had a heart attack, and died.  His
{\it Opere matematiche\/} in three volumes (Hoepli, Milano 1914, 1915,
and 1917) were edited under Bertini's direction.

In 1876, in Pisa, Bertini classified the plane Cremona transformations
that are equal to their own inverses, the plane involutions of order
two.  He proved that they decompose into products of irreducible
involutions, each of which is equivalent under a Cremona transformation
to one of only four types, a surprisingly simple result.  The first two
types were already well known: reflections in lines, and de Jonqui\`eres
involutions.  However, Bertini showed that each de Jonqui\`eres
involution is irreducible, and is given by a curve of degree $n$, for
some $n$, with a multiple point of order $n-2$; a given point, in
general position in the plane, corresponds to its harmonic conjugate
with respect to two other points, namely, the two points of the curve on
the line determined by the given point and the multiple point.  An
involution of the third type is given by a net of cubics through seven
points; the cubics through an eighth point in general position obviously
pass as well through a ninth, and the ninth corresponds to the
eighth.  This type had been found ten years earlier by Geiser, but
described differently.

An involution of the fourth type is given by a 3-dimensional system of
sextics passing doubly through eight points; the sextics through a ninth
point in general position can be proved to pass as well through a tenth,
and the tenth corresponds to the ninth.  These involutions were new, and
have become known as {\it Bertini involutions}.  Later, in 1889, Bertini
gave a simple geometric proof that an involution of each of the four
types is rational; that is, the pairs of corresponding points form the
fibers of a rational 2-to-1 map onto the plane.  (This result, in
essence, had already been stated by Noether in 1876, and was also proved
by L\"uroth in 1889, but algebraically.  In 1893, Castelnuovo proved the
more general result that any rational transform of the plane is
birational to the plane.)

Bertini's classification represented a philosophical break with Cremona,
who saw his transformations only as a tool for reducing the complexity
of given geometric figures, and not as objects of study in their own
right.  Cremona received Bertini's work coldly, unjustifiably so.
Bertini explained it to him in person before the appearance of his main
paper \cite{Bert77}.  However, Cremona said only that he had already
discovered the fourth type of involution, and had communicated its
existence in a letter to his former student, Ettore Caporali (1855--86).
Bertini respectfully added Cremona's existence proof to the page proofs;
see \cite{Bert77, p.~273}.  Moreover, Bertini needed to make a technical
restriction on the fundamental locus, and did so openly in the first
paragraph of his paper \cite{Bert77}.  The restriction was eventually
removed by Castelnuovo and Enriques.

Meanwhile, two years later, Caporali took up the classification from a
different point of view, and suggested that there might be, in addition
to Bertini's four types, infinitely many others.  The significance of
this entire incident is suggested by the length and strength of the
discussions of it by Castelnuovo \cite{Cast33, pp.~746--7} and by Scorza
\cite{Scor34, pp.~105--7} nearly fifty years after the fact.  In
particular, Castelnuovo called Caporali's point of view ``much less
interesting" than Bertini's, which ``opened new horizons to algebraic
geometry.''  And Scorza concluded with this bit of wisdom: it takes
longer to appreciate a conceptual advance when it is more original and
more profound.

Bertini's research reflects, by and large, Cremona's influence in its
subject and its style.  Both geometers employed synthetic methods and
analytic methods with equal facility.  Both wrote succinct, precise,
elegant treatments, which reflect their powerful, penetrating intellects.
Many of their investigations concerned special properties of given
figures; these works have lost much of their interest, but had their
importance in pointing the way to the general theory.

Much of Bertini's research grew out of his lecture preparations.  For
him, research and teaching were two aspects of the same activity.  He
kept current by reworking the latest advances, often putting them in a
new and perspicacious form, which he included in his lectures and
papers.  Here are some examples.

In 1881 Bertini enumerated the 5-secant conics to the quintic space
curve.  In 1884 he wrote his extensive and elaborate memoir on the
cubic surface, with its 27 lines and 45 tritangent planes.  In 1888 he
gave a simple and ingenious geometric proof of the following one of
Noether's theorems: any plane curve can be transformed into another one
that has only ordinary multiple points, via Cremona transformations of
the ambient plane.  As a consequence, he extended Pl\"ucker's formulas
to plane curves with arbitrary singularities.  Later, in 1891 he
transformed a given curve, via arbitrary birational transformations,
into a plane curve with ordinary double points; this result had already
been in use for a long time, but not yet proved rigorously.  In 1890
and in 1908, he advanced the theory of linear series on an abstract
curve.  In 1896 he studied the 21 quadruple tangents of the Cayleyian
of a plane quartic.  In 1898 he studied pencils of quadrics and the
linear spaces on a quadric of even dimension.

In 1880 Bertini left Pisa for Pavia, exchanging positions with Riccardo
De Paolis (1854--92), another one of Cremona's former students.  De
Paolis wanted to be in Pisa, and Bertini wanted to be with his two close
colleagues, Eugenio Beltrami (1835--1899) and Casorati.  Beltrami had
taught algebra and analytic geometry at Bolognia when Bertini was a
student there.  In 1890 Casorati died, and Bertini took over the
responsibility for teaching analysis.  Two years later, De Paolis died
prematurely, and Bertini reclaimed his old position in Pisa on the
encouragement of Luigi Bianchi (1856--1928) and of Dini.  Bertini
remained in Pisa for the rest of his life.  Among his students during
this period were Carlo Rosati (1876--1929), Ruggiero Torelli
(1884--1915), and Giacomo Albanese (1890--1947).

Bertini was a scrupulous and zealous teacher, who viewed teaching as a
ministry.  Each year he covered a different subject.  He prepared his
lectures with great care, both in plan and in detail, obtaining a
precise, efficient, clean development.  His delivery was clear and dry,
yet lively.  He often asked questions and proposed problems to engage
his audience, and he knew how to maintain decorum.  He demanded hard
work of his students, but they loved him, and called him ``Pap\'a
Bertini.''

Bertini was a tall heavyset man.  He had a fresh open mind, and a
noble upright character.  He was modest in every way, and found a
friend in every one.  He was kind and affectionate.  He was frank,
optimistic, and generous.  He was full of advice, encouragement, and
help.  He had a strong sense of duty, and an inflexible sense of
justice.  He was looked up to as an inspiring example of high morality.

Bertini's lectures on projective algebraic geometry at the Ateneo pisano
were written up and lithographed during the academic year 1889--99 by
his assistant and former student, Scorza.  They were then revised,
expanded with an appendix, and published in 1907 as the book
\cite{BertI}.  It was republished in 1923, and translated into German in
1924.  The book contains all the essential results about projective
varieties that had been obtained in the preceding decades.  The appendix
treats algebraic curves and their singularities; it was based in part on
the extensive summaries prepared by Corrado Segre (1863--1924) for his
own courses.  Bertini's exposition is systematic, extensive, and lucid.
It was the crowning achievement of his scientific work, and is still
studied and cited today.

On 1 August 1922 Bertini had to retire because of age (75).  He was
succeeded by Rosati, who had been his assistant.  However, for the next
ten years, Bertini continued to teach for free, as professor emeritus,
because teaching was so important to him.  He also turned these lectures
into a text book \cite{BertC}, which he wrote himself at the age of 82.
The book was aimed at the second and third year students; it discusses
assorted topics in the geometry of projective space, such as quadratic
transformations of the plane, the formulas of Pl\"ucker and Cayley
relating numerical characters of plane and space curves, various types
of line complexes, the cubic surface, and the Steiner surface.  The book
nobly closed Bertini's scientific production.  He gave his scientific
library to the faculty of science, and funded an annual prize for the
best laurea dissertation in pure mathematics.

On 24 February 1933, Eugenio Bertini died in Pisa after a brief illness,
survived only by his daughter, Eugenia.  He also had a son, Giulio, with
his wife Giulia, n\'ee Boschi.  She died in Pisa on 23 January 1915
after a long and difficult illness, and Giulio died tragically, 27 years
old, on 22 September 1922 in a traffic accident.  Bertini was buried, as
he wished, without the pomp and ceremony of a traditional large outdoor
funeral, in the cemetery of Forli, along side his wife.  Nevertheless,
a year later, on 24 February 1934, a remarkable number of people,
including immediate family, former students, close colleagues, academic
administrators, and political dignitaries, met at the R. Istituto
Tecnico di Forli to celebrate Bertini's life and to unveil a memorial
plaque; the event is lovingly described by Conti in his chronicle
\cite{Con34}.  Throughout his life, Bertini dedicated himself humbly and
altruistically to family, school, and science.

\sct The original first theorem

In this section and the next, we'll discuss Bertini's first fundamental
theorem on linear systems, which concerns the multiple points on the
general member.  Originally, in 1882, Bertini \cite{Bert82} worked only
with systems on the full $n$-dimensional projective space $\I P^n$ over
the complex numbers, and this is the only case that we'll consider case
here.  In the next section, we'll go on to discuss the extension of the
theorem to an arbitrary ambient variety in arbitrary characteristic.
In this section, we'll begin by reviewing some basic notions in the
form in which they were considered in Bertini's time.  Then we'll
consider the content of the theorem.  Finally, we'll go through
Bertini's proof, which is simple, rigorous, and interesting.

On the projective space $\I P^n$, a {\it linear system\/} is simply the
family of all hypersurfaces $U_{\I t}$ of the form,
        $$U_{\I t}:t_0u^{(0)}+\cdots+ t_su^{(s)}=0,$$
 where $\I t:=(t_0,\dots,t_s)$ is a point of a projective space $\I
P^s$, the {\it parameter space}, and where $u^{(0)},\dots,u^{(s)}$ are
linearly independent homogeneous polynomials of the same degree in
$n+1$ variables $x_0,\dots,x_n$.  The system is called a {\it (linear)
pencil\/} if $s=1$ and a {\it net\/} if $s=2$.

 If all the members $U_{\I t}$ contain a common hypersurface $U:u=0$,
then $U$ is said to be {\it fixed}, and the system defined by the
quotients $u^{(0)}/u,\dots,u^{(s)}/u$ is called the {\it residual}
system.  There is a smallest residual system: it is the one where $u$ is
the greatest common divisor of the $u^{(i)}$.  If all the members of
this smallest residual system contain a common point or a common
variety, then it is called a {\it base point\/} or {\it base variety\/}
of the original system.  So a base variety has dimension at most $n-2$.
(Nowadays, it is more common to consider a fixed component to be a base
variety as well.)

Let $U:u=0$ be a hypersurface.  A point of $U$ at which all the partial
derivatives of order $r-1$,
        $$u_{i_1\dots i_{r-1}}:=
         {\partial^{r-1}u\over\partial x_{i_1}\dots\partial x_{i_{r-1}}}
                        \hbox{ where }0\le i_j\le n,$$
 vanish is called a point of {\it multiplicity\/} $r$, or an {\it
$r$-fold point}.  It may also be an $(r+1)$-fold point (although some
authors, including Bertini, sometimes exclude this possibility as part of
the definition).  It is called a {\it multiple\/} point, or a {\it
singular\/} point, if $r\ge2$, but the exact value of $r$ is
unimportant.  If a point is not multiple, then it is called {\it
simple\/}.  Given $r$, the various $r$-fold points of the various
members $U_{\I t}$ of the linear system form a closed subset $M_r$ of
$\I P^n\x\I P^s$, and a given $r$-fold point is said to be {\it
variable} if it lies in an irreducible closed subset of $M_r$ that
covers $\I P^s$.

The preceding concepts are illustrated in the two simple examples shown
in \fig1.  (All the figures were drawn using the \TeX-PostScript macro
package {\tt PSpictures} written by Thorup.)  In both examples, the
ambient space is the plane $\I P^2$ with inhomogeneous coordinates
$x,y$, and the parameter space is the line $\I P^1$ with inhomogeneous
coordinate $t$.  In the first example (the one on the left), each curve
$U_{\I t}$ consists of two components: one is the vertical line through
the point $\I x(t,0)$, and the other is the $x$-axis with multiplicity
two.  The latter is a fixed component; $\I x$ is a variable $3$-fold
point, or triple point; and the point at infinity on each vertical line
is a base point, and the only one.  In the second example, the vertical
lines of the first are replaced by the pencil of lines through the
origin.  So it is a base point, and the only one; moreover, it is a
variable triple point, although $\I x$ is fixed in the plane.  Once
again, the $x$-axis is a fixed component with multiplicity two.

 \PSunit=1pc \newdimen\Xdim
 \def\SetDash{\setdash{2.4\PSlinewidth}}
\centerline{\eightpoint
 \def\fixpart{
  \PSarrow{-3.5 0 4.75 0} \TX{ $x$\hss}{4.75 0}
  \PSarrow{0 -1.75 0 2} \TX{ $y$\hss}{0 1.75}
  \setPSlinewidth{1pt}
  \PSline{-3 0.1 3.75 0.1}  \PSline{-3 -0.1 3.75 -0.1}
 \TX{\box0\hss}{0.75 0.75}
 }
 \setbox0=\hbox{ $U_{\bf t}:(x-t)y^2=0$}
 \Xdim=\wd0 \advance\Xdim3.75\PSunit
 \PSpicture{\Xdim}{4.25\PSunit}\llcoords{-3 -1.75}
   \fixpart
  \PSline{0.75 -1.5 0.75 1.5}
  \setPSlinewidth{2pt}\PSmark{0.75 0}\setPSlinewidth{1pt}
  \TX{ $\I x(t,0)$\hss}{0.75 -0.75}
  \setPSlinewidth{2pt}\PSmark{-2 0}\setPSlinewidth{1pt}
  \SetDash \PSline{-2 -1.5 -2 1.5} \nodashes
 \endPSpicture
 \hfil
 \setbox0=\hbox{ $U_{\bf t}:(tx-y)y^2=0$}
 \Xdim=\wd0 \advance\Xdim3.75\PSunit
 \PSpicture{\Xdim}{4.25\PSunit}\llcoords{-3 -1.75}
  \fixpart\PSline{-1.5 -1.75 1.5 1.75}\TX{ $\I x(0,0)$\hss}{0.75 -0.75}
  \SetDash\PSline{1.5 -1.75 -1.5 1.75} \nodashes
  \setPSlinewidth{2pt}\PSmark{0 0}
 \endPSpicture
}\par\nobreak
 \centerline{\eightpoint\fig1.
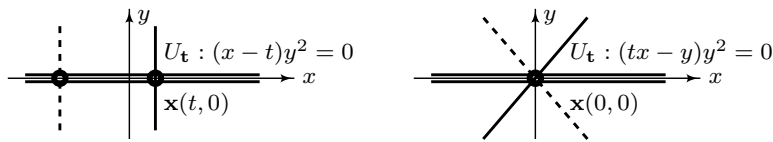
  Fixed components, base points, and variable triple points}
\medbreak

In general, when $n\ge2$, each $U_{\I t}$ will have an $r$-fold point
where $r\ge2$ if there is a hypersurface $U:u=0$ such that
$(r-1)U:u^{r-1}=0$ is fixed.  Indeed, then $U_{\I t}$ is defined by an
equation of the form $u^{r-1}v_{\I t}=0$, and all the partial
derivatives of order $r-1$ of the product $u^{r-1}v_{\I t}$ vanish at
every point of the intersection $W_{\I t}$ of $U$ and $V_{\I t}$, where
$V_{\I t}:v_{\I t}=0$ is the residual hypersurface.  So the points of
$W_{\I t}$ are $r$-fold points of $U_{\I t}$, and each component $W'_{\I
t}$ of $W_{\I t}$ has dimension at least $n-2$.  Suppose that some
$W'_{\I t}$ is not contained in the base locus, and let $U'$ be any
component of $U$ that contains $W'_{\I t}$, but not $V_{\I t}$.  Then
each point $\I x$ of $U'$ must lie on some member $V_{\I t'}$ of the
residual system (otherwise $W'_{\I t}\x\I P^s$ would be a component of
the intersection of $U'\x\I P^s$ with the total space of the residual
system); so $\I x$ is an $r$-fold point of $U_{\I t'}$.  Hence $U'$
consists entirely of variable $r$-fold points.

This is not the only way that there can be a variable $r$-fold point
off the base locus, although Bertini asserted that it was so in his
original formulation of his first theorem \cite{Bert82, \p.26}.  For
instance, if there is a fixed hypersurface with an $r$-fold point, off
the base locus or not, then this point will be a variable $r$-fold.  In
his book \cite{BertI, \p.227}, Bertini was more careful, and put the
theorem essentially as follows.
        \proclaim Theorem 1 \(On variable $r$-fold points).  On $\I
P^n$, a variable $r$-fold point of a linear system is an $(r-1)$-fold
point of every member.
        \endproclaim

The most important case of Theorem~(3.1) occurs when $r$ is $2$ and
there are no fixed components.  Then the result is always put in
contrapositive form; Bertini himself \cite{Bert82, \p.26} did so
essentially as follows.
        \proclaim Theorem 2 \(On variable singular points).  On $\I
P^n$, if a linear system has no fixed components, then a general member
has no singular points outside the base locus.
        \endproclaim

By a {\it general\/} member of a linear system is meant one that is
represented by a point in the parameter space ${\I P^s}$ that satisfies
no special algebraic-geometric conditions; in other words, the point
lies in no subvariety that depends on the issue in question.  Curiously,
Bertini changed terminology over the years, doubtless reflecting the
prevailing usage.  Thus, in his 1877 paper \cite{Bert77, \p.246}, he
wrote ``curva generale''; in his 1882 paper \cite{Bert82, \p.26}, he
wrote ``spazio arbitrario''; and in his book \cite{BertI, \p.227}, he
wrote ``ipersuperficie generica.''

Without the hypothesis of linearity of the system, the theorem may
fail.  A simple counterexample is given in \fig2.  Here the curves
$U_{\bf t}$ are the horizontal translates of an irreducible plane cubic
with a cusp at the origin.  Thus each curve has a point of multiplicity
two, or double point, but the various curves have no point in common,
not even at infinity.

\centerline{
  \setbox0=\hbox{ $U_{\bf t}:y^2=(x-t)^3$}
 \Xdim=\wd0 \advance\Xdim4.5\PSunit
 \PSpicture{\Xdim}{4\PSunit}\llcoords{-1 -1.5}
  \PSarrow{-1 0 6 0} \TX{ $x$\hss}{6 0}
  \PSarrow{0 -1.75 0 2} \TX{ $y$\hss}{0 1.75}
  \setPSlinewidth{1pt}
  \def\cubic{\PS{newpath 0 0 moveto 0 0 1 0 1.5 1.5 curveto stroke}
        \PS{newpath 0 0 moveto 0 0 1 0 1.5 -1.5 curveto stroke}
        \PSmark{0 0}
  }
  \SetDash \cubic \nodashes
  \llcoords{-3 -1.5} \cubic  \TX{\box0\hss}{1.5 1}
 \endPSpicture
}\par\nobreak
 \centerline{\eightpoint
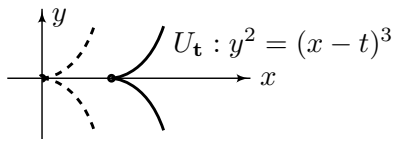
 \fig2. A nonlinear system with a variable double point}
 \medbreak

Theorem~(3.2) had already been proved in 1871, as Bertini \cite{Bert82,
p.26} noted, by Rosanes \cite{Ros71, \p.100} in a particular case, that
of a {\it homoloidal net} of plane curves.  This is the case where
$s=2$, and the $u^{(i)}$ have no common factor.  Moreover, they must
define a Cremona transformation of the plane; that is, the rational map,
 $$\I x\mapsto \big(u^{(0)}(\I x), u^{(1)}(\I x), u^{(2)}(\I x)\big),$$
 is a birational transformation of the plane onto itself.  Rosanes'
computational proof is special to this case, and gives no hint of a
more general result.

  In 1877 Bertini stated Theorem~(3.2) for an arbitrary linear system
of plane curves on \p.246 of his main paper \cite{Bert77} on plane
involutions.  It may be true, as Castelnuovo \cite{Cast33, \p.747} put
it, that Bertini established this result ``by means of intuitive
considerations.''  However, these considerations do belong to the
standard theory of envelopes, which can be justified using calculus, as
Bertini did in his book \cite{BertI, Fn.~1, \p.225}.  Nevertheless, in
the book, Bertini gave a somewhat different proof of the theorem: he
repeated the proof in his paper \cite{Bert82}, see below.  On the other hand,
such complete analytic local proofs of the theorem in the plane are
given in the texts of Picard and Simart \cite{PiSi06, \p.51}, of
Francesco Severi (1879--1961) \cite{Sev08, \p.27} and of Enriques and
Chisini \cite{EnCh18, \p.181}.  In fact, Severi proved the full
Theorem~(3.1) itself in the plane, although he attributed only the
special case of Theorem~(3.2) to Bertini.  A more general version of
this local proof is given at the end of the next section.

Theorem~(3.1) applies equally well to a variable $r$-fold base point
with $r\ge3$.  A simple example with $r=3$ is shown in \fig3.  Here
$U_{\I t}$ is the cone in $\I P^3$ with vertex at $\I x_{\I t}:=(0,0,t)$
over a fixed irreducible curve $C$ in the $(x,y)$-plane with a double
point at the origin.  The vertex $\I x_{\I t}$ is a triple point of
$U_{\I t}$, but a double point of $U_{\I t'}$ for $\I t'\neq\I t$.  As
$\I t$ varies, $\I x_{\I t}$ sweeps out the $z$-axis, which is the full
base locus.
        \medskip
\centerline{\PSunit=2pc
 \PSpicture{4.5\PSunit}{5.75\PSunit}\llcoords{-1 -2.25}
  \PSarrow{-1 0 3.5 0} \TX{ $x$\hss}{3.5 0}
  \PSarrow{-1.5 -2.25 1.5 2.25} \TX{ $y$\hss}{1.5 2.2}
  \PSarrow{0 -2.25 0 3.5} \TX{ $z$\hss}{0 3.35}
  \TX{\hss$\I x_{\I t}:=(0,0,t)$ }{0 2.75}
  {\setPSlinewidth{1pt}
   \PS{newpath 0 0 moveto 0.0 0 2 -0.24 4.5 2.25 curveto stroke}
   \PSline{0 2.75 4.5 2.25}
   \PS{newpath 0 0 moveto 0.0 0 2.36 0.24 1.5 -2.25 curveto stroke}
   \PSline{0 2.75 1.5 -2.25}
   \PSline{0 2.75 0 0}
  \PSmark{0 2.75}
   \SetDash
    \PSline{0 1.5 4.5 2.25}
    \PSline{0 1.5 1.5 -2.25}
        {\setPSlinewidth{1.5pt}\PSline{0 1.5 0 0}}
   \nodashes
  \PSmark{0 1.5}
  }
  \TX{$C:y^2=x^3$\hss}{2 -1.125}
  \TX{ $U_{\bf t}:y^2(z-t)=tx^3$\hss}{0.5 3}
 \endPSpicture
}\par\nobreak
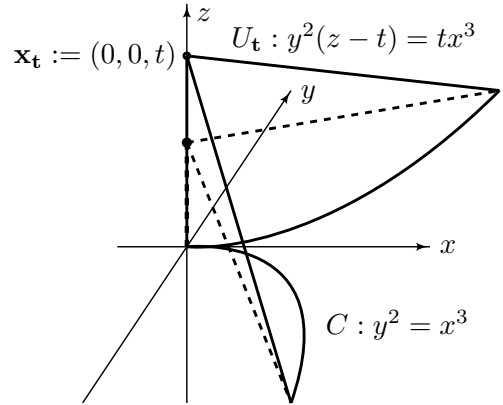
 \centerline{\eightpoint\fig3. A linear system with a variable triple
point in the base locus}
        \medbreak Bertini proved his theorem rigorously, more or less as
follows, although he did not pay as much  attention to the details in
the first part.  Let $\I x\in U_{\I t_0}$ be the variable $r$-fold
point; by definition, $(\I x,\I t_0)$ lies in a component $M$ of $M_r$
such that the projection $\mu\:M\to\I P^s$ is surjective.  Hence $\mu$
is smooth on some dense open subset $M^0$ because the characteristic is
zero.  Now, we have to prove that $\I x$ lies in the closed subset $N$
of $\I P^n$ where all the partial derivatives of order $r-2$ of
$u^{(k)}$ vanish for $0\le k\le s$.  Since $M^0$ is dense in $M$, it is
enough to prove that the image of $M^0$ is contained in $N$.  Thus we
may assume that $\mu$ is smooth at $\I x$.

We may replace $u^{(0)},\dots,u^{(s)}$ by any invertible linear
combination of themselves; so we may assume that $U_{\I t_0}:u^{(0)}=0$.
Given $k$, set
        $$v(t):=u^{(0)}+tu^{(k)}.$$
  Since $\mu$ is smooth at $\I x$, there are power series
$x_0(t),\dots,x_n(t)$ such that
  $$\displaylines{\I x(t):=(x_0(t),\dots,x_n(t))\in M,\ \I x=\I x(0),\cr
                \mu(\I x(t))=(1,0,\dots,0,t,0,\dots,0)\cr}$$
 where the last $t$ is placed in the $k$th coordinate.  Then all the partial
derivatives of order $r-1$, with respect to the $x_i$, of $v(t)$ vanish at
$\I x(t)$; in other words,
        $$u^{(0)}_{i_1\dots i_{r-1}}(\I x(t))
        +tu^{(k)}_{i_1\dots i_{r-1}}(\I x(t))=0
        \hbox{ where }0\le i_j\le n.\eqno(*)$$
 Equation $(*)$ holds identically in $t$.  So differentiating it yields
        $$\displaylines{\qquad
                \sum_lu^{(0)}_{i_1\dots i_{r-1}l}(\I x(t)){dx_l(t)\over dt}
  +t\sum_lu^{(k)}_{i_1\dots i_{r-1}l}(\I x(t)){dx_l(t)\over dt}\hfill\cr
  \hfill      +u^{(k)}_{i_1\dots i_{r-1}}(\I x(t))=0.\qquad\cr}$$
 Multiply by $x_{i_1}(t)$ for $i_1=0,\dots,n$, and take the sum.  Let
$m$ be the common degree of $u^{(0)}$ and $u^{(k)}$.  Then Euler's
identity yields
        $$\displaylines{\qquad
                (m-r+1)\sum_l\Bigl(u^{(0)}_{i_2\dots i_{r-1}l}(\I x(t))
        +tu^{(k)}_{i_2\dots i_{r-1}l}(\I x(t))\Bigr){dx_l(t)\over dt}\hfill\cr
        \hfill+(m-r+2)u^{(k)}_{i_2\dots i_{r-1}}(\I x(t))=0.\qquad\cr}$$
 However, the expression between the large parentheses vanishes by
$(*)$.  Hence $u^{(k)}_{i_2\dots i_{r-1}}(\I x(t))=0$.  Setting $t=0$
yields $u^{(k)}_{i_2\dots i_{r-1}}(\I x)=0$, which was to be proved, and the
proof of Theorem~(3.1) is complete.

\sct The extended first theorem

In this section, we'll discuss the extension of Bertini's first theorem
to an arbitrary ambient variety.  We'll trace the evolution of its
statement and of its proof in the hands of Bertini and others.  In
particular, we'll discuss some rather sketchy and intuitive proofs given
by Enriques and by Severi.  These proofs will only be described briefly
to give their flavor.  Doubtless, with some effort, they could be
completed and made rigorous, but doing so would provide involved proofs
that offer no new insight.  Until explicit mention is made to the
contrary about half way through the section, the ground field will be
the complex numbers.

To date, Theorem~(3.1) has been extended only in the special case of
Theorem~(3.2), where the minimal order of multiplicity $r$ is 2.
However, at the end of this section, by refining, developing and
completing some of the old arguments, we'll prove a full extension,
Theorem~(4.2), of Theorem~(3.1) itself.

In the extensions of Theorem~(3.1) and Theorem~(3.2), the original
statements must be refined to make allowance for the possible presence
of singularities in the ambient variety, which is arbitrary now.
Thus Theorem~(3.2) acquires the following familiar statement.
        \proclaim Theorem 1 \(On variable singular points, extended).
On an arbitrary ambient variety, if a linear system has no fixed
components, then the general member has no singular points outside of
the base locus of the system and of the singular locus of the ambient
variety.
        \endproclaim

Here the ambient variety is abstract, and need not be embeddable in any
projective space, although only embeddable varieties were considered
before about 1950.  However, this matter is of no real importance since
the question is local.

The ambient variety was a surface when Theorem~(4.1) was stated for the
first time and by Enriques.  He did so in 1893 in his first major work
\cite{En93, \p.42}, which initiated a general theory of linear systems
(notably including adjoint systems) on an abstract surface.  For
Enriques, a linear system was obtained by taking a suitable projective
model of the surface, cutting it by the members of a linear system of
hypersurfaces, and then stripping away some, or all, of the {\it fixed
components}, which are the components common to all the members of the
induced system.  Thus the members of a linear system are {\it divisorial
cycles}, or linear combinations of subvarieties of codimension 1.  These
cycles need not be defined by a single equation locally about any given
singular point of the ambient surface.  However, this matter too is of
no real importance since the theorem asserts nothing about these
singular points.

When Enriques stated the theorem, he did not simply assume that the
system has no fixed components; curiously, he also assumed that a
general member is irreducible.  Enriques did not dwell on the details of
his proof, but in brief he argued by contradiction as follows (and
repeated the argument, in slightly more detail, in his 1896 paper
\cite{En96, \pp.230--31}): if the theorem were false, then the given
system would contain a pencil with a variable singular point, and the
surface would have a birationally equivalent model in $\I P^3$ so that
the equivalent pencil consists of the sections by the planes through a
line; however, each of these planes would be tangent to the model, and
this situation is absurd.

The theorem was extended further, to an ambient variety $X$ of arbitrary
dimension, by Severi in his 1906 paper \cite{Sev06, \p.169}.  On \p.168,
he credited Bertini's 1882 paper \cite{Bert82} for the original
statement and Enriques 1896 paper \cite{En96, \pp.230--31} for the
extension to an arbitrary surface.  At the same time, Severi
acknowledged the receipt of a letter from Bertini, in which Enriques'
extension was used.  Severi then sketched a proof, which is in some ways
similar to and in some ways different from Enriques' proof and also from
Bertini's proof in his book \cite{BertI, \S18, \p.239}, which had
already appeared in lithograph in 1899, but not yet in print.
Generously, Severi \cite{Sev06, Oss., \p.170} described his proof as
nothing but a geometric form of Bertini's 1882 proof.

Severi's sketch looks much like this.  Let $\I x$ be a variable singular
point at which $X$ is smooth.  We may assume that $X$ is embedded in a
$\I P^n$ as a hypersurface and that the linear system on $X$ is cut out,
off a fixed component not containing $\I x$, by a linear system of
hypersurfaces $U_{\I t}$.  If this system has $\I x$ as a variable
singular point too, then it is a base point by Bertini's Theorem~(3.2)
for the ambient $\I P^n$.  Otherwise, a general $U_{\I t}$ is smooth
at $\I x_{\I t}$; replace it by its tangent hyperplane, use the dual
variety to conclude that this hyperplane contains $\I x$, and finally
deduce that every $U_{\I t}$ too contains $\I x$; thus again $\I x$ is a
base point.

As noted above, Bertini himself had already proved the extended theorem
in arbitrary dimension in his book \cite{BertI}.  Curiously, Bertini did
not cite Enriques or Severi, although he did cite many other relevant
works at appropriate places through out the book.  Bertini's proof and
Severi's are basically alike.  However, Bertini did not require $X$ to be
a hypersurface in the ambient $\I P^n$.  On the other hand, he too used
his Theorem~(3.2) for $\I P^n$ to reduce to the case that the general
$U_{\I t}$ is smooth at $\I x_{\I t}$.  However, he did not then replace
this $U_{\I t}$ by its tangent hyperplane; rather, he gave a direct
analytic argument that every $U_{\I t}$ must contain $\I x$.  In fact,
it is unnecessary to make this application of his original theorem for
either of the following two reasons: first, if we re-embed $X$ using a
Veronese embedding, then we may assume that all the $U_{\I t}$ are
hyperplanes, so smooth; second, Bertini's analytic argument works even
when the general $U_{\I t}$ is not smooth, as we'll see now.

Bertini's analytic argument is rigorous here again, and so it provides a
rigorous proof of Theorem~(4.1).  The argument is much like that in his
proof of Theorem~(3.1), but there is a twist.  The argument runs as
follows.  The setup here is the same as that in the proof of
Theorem~(3.1), except that now the ambient variety $X$ is a subvariety
of $ \I P^n$ and the parameterized variable point $\I x(t)$ lies outside
the smooth locus of $X$ and inside the singular locus of $X\cap U_{\I
t}$, where $U_{\I t}:v(t)=0$ is a parameterized hypersurface in $\I
P^n$.  These conditions imply that the vector,
        $$\left({dx_0\over dt}(0),\dots,{dx_r\over dt}(0)\right),$$
 lies in the tangent space to $X$ at $\I x=\I x(0)$; they also imply
that this tangent space lies in the tangent hyperplane to $U_{\I t_0}$ at
$\I x$ if $U_{\I t_0}$ is smooth there.  So, in any event,
        $$\sum_iv_i(0)(\I x){dx_i\over dt}(0)=0
                \hbox{ where }v_i:=\pd v{x_i}.$$
 On the other hand, $v(t)(\I x(t))=0$.  Differentiating this equation
yields
        $$\sum_iv_i(0)(\I x){dx_i\over dt}(0)+\pd vt(\I x)=0.$$
 Hence $\pd vt(\I x)=0$.  Now $v(t):=u^{(0)}+tu^{(k)}$.  So $u^{(k)}(\I
x)=0$, as was required.  Thus the proof of Theorem~(4.1) is complete.

Zariski gave an intrinsic proof of Theorem~(4.1) in his paper
\cite{Zar44} of 1944.  (Zariski was in Italy from 1921 to 1927, and
studied with Castelnuovo, Enriques and Severi; for a lovely description
of Zariski's stay, see \cite{Par91}.)  Ten years earlier, Enriques had
given an intrinsic proof for surfaces in his book with Campedelli
\cite{EnCa32, \p.31} (again assuming a general member is irreducible);
Enriques simply asserted that the theorem can be proved in essentially
the same way as it was for the plane in his book with Chisini
\cite{EnCh18, \p.181} (where a general member is not assumed to be
irreducible), because this proof is local in nature.  Enriques'
assertion was repeated by Zariski in his book \cite{Zar35, \p.26} of
1935.

Zariski's 1944 proof is conceptually more advanced, and runs as follows.
Replacing the ambient space $X$ by an open subset, we may assume that
$X$ is smooth and that the system has no fixed components and no base
points.  Then the total space of the system is smooth too, because it is
a locally trivial bundle of projective spaces; indeed, the fiber over a
point of $X$ is just the set of members that contain this point.  A
general member of the system is, therefore, smooth everywhere by what is
now usually called ``Sard's theorem,'' which holds in characteristic 0.

Zariski proved Sard's theorem in the case at hand; in fact, the proof
works in general.  Namely, Zariski noted that the ``generic'' member of the
system is regular (an intrinsic property, defined by means of
uniformizing parameters), because, as a scheme, it is the fiber over the
``generic'' point of the parameter space $\I P^s$; so its local rings are
simply localizations of the local rings of the total space.  Now, in
characteristic zero, a regular scheme is smooth \cite{Zar47, Thm.~7},
and smoothness is an open condition, as it is defined ``by means of
nonvanishing Jacobians'' \cite{Zar44, \p.139}.

It was Zariski's work on Bertini's two theorems that helped motivate
him to study the relationship between regularity and smoothness,
according to Mumford \cite{Mum72, \p.5}.  In turn, it is clear from
what Zariski himself wrote, in his paper \cite{Zar44r, \p.474} on the
reduction of singularities of 3-folds, that this work on singularities
had motivated him to work on Bertini's first theorem.

Zariski \cite{Zar44, \p.140} also provided (essentially) the following
two examples, which show that Theorem~(3.2), as stated, is false in
positive characteristic $p$.  (Earlier, van der Waerden \cite{vdW35,
\S3., \p.133} had noted that his version of Enriques' 1896 proof
required an appropriate hypothesis of separability.)  In both examples,
the ambient space is the plane $\I P^2$ with inhomogeneous coordinates
$x,y$, and the parameter space is the line $\I P^1$ with inhomogeneous
coordinate $t$.  In the first example, the system consists of the
curves,
        $$U_{\I t}:x^p-ty^p=(x-t^{1/p}y)^p=0,$$
 which consist entirely of $p$-fold points.  In the second example,
$p\neq2$, and the system consists of the (absolutely) irreducible
curves,
        $$U_{\I t}:y^2-x^p+t=y^2-(x-t^{1/p})^p=0,$$
 which have a variable double point at $(0,t^{1/p})$.

Finally, here is the new extension of Bertini's Theorem~(3.1), which
treats variable $r$-fold points for an arbitrary $r$.  This extension,
Theorem~(4.4) below, is valid for an arbitrary ambient variety $X$ in
arbitrary characteristic $p$.  Validity for $p>0$ is made possible by a
suitable new notion, ``separable variable $r$-fold point,'' which is
defined below.  Of course, the separability requirement is automatically
satisfied when $p=0$.

The various $r$-fold points of the members of a linear system on $X$,
parameterized by $\I P^s$, form a subset, $M_r$ say, of $X\x\I P^s$.  A
given $r$-fold point    will be called  {\it variable}, resp.,
{\it separable variable}, if it lies in a irreducible subset $M$ of
$M_r$ that is closed in $X\x\I P^s$ and that projects onto $\I P^s$, resp.,
and that projects separably onto $\I P^s$.  In fact, the proof of
Theorem~(4.2) shows that $M_r$ is closed in $X\x\I P^s$; so $M$ is
closed in $M_r$.  Notice that the multiple points in Zariski's two
examples are variable, but not separable.
        \proclaim Theorem 2 \(On variable $r$-fold points, extended).
Let $X$ be a variety over an algebraically closed field of arbitrary
characteristic, and $\I x$ a point of $X$ outside of its singular locus.
If $\I x$ is a separable variable $r$-fold point of a linear system, then
$\I x$ is an $(r-1)$-fold point of every member.
        \endproclaim

Since the matter is local on $X$, we may replace $X$ by an affine
neighborhood of $\I x$ on which there exist ``uniformizing coordinates'';
these are functions $\xi_1,\dots\xi_n$, where $n:=\dim X$, whose
differentials $d\xi_i$ form a free basis of the 1-forms.  Identify
$\xi_i$ with its pullback to $X\x X$ via the second projection, denote
by $\eta_i$ its pullback via the first projection, and set
$\delta_i:=\eta_i-\xi_i$.  Given any (regular) function $u$ on $X$,
denote by $u(\eta)$ its pullback to $X\x X$ via the first projection,
and by $u(\xi)$ that via the second.  Then $u(\eta)$ has a unique
expansion,
 $$u(\eta)=\sum_{\scriptstyle0\le i_j\le n\atop\scriptstyle0\le q<r}
        u_{i_1\dots i_q}(\xi)\delta_{i_1}\cdots\delta_{i_q}+w,$$
 where the $u_{i_1\dots i_q}$ are suitable functions on $X$ and where
$w$ is a suitable function on $X\x X$ that lies in the $r$th power of
the ideal $\Id$ of the diagonal.  In characteristic 0, the $u_{i_1\dots
i_q}$ are simply scalar multiples of the partial derivatives of $u$,
and this expansion of $u$ is essentially the Taylor expansion.  In
arbitrary characteristic, the expansion can be obtained by plugging in
$\delta_i+\xi_i$ for $\eta_i$ in $u(\eta)$ and collecting the terms in
the products $\delta_{i_1}\cdots\delta_{i_q}$; the expansion is unique
because these products provide a free basis of $\Id^q/\Id^{q+1}$ as a
module over the ring of functions on $X$ because $\Id/\Id^2$ is equal
to the module of 1-forms.

A point $\I x$ of $X$ is an $r$-fold point of the divisor $U:u=0$ if and
only if all the functions $u_{i_1\dots i_q}$ vanish at $\I x$.  In
characteristic 0, this statement is essentially the definition in
Section~3 for the case $X=\I P^n$.  In general, we must see that this
vanishing holds if and only if $u\in \Mx^r$ where $\Mx$ is the
maximal ideal of $\I x$.  To do so, identify $X$ with $X\x\{\I x\}$
inside $X\x X$.  Then $u(\eta)$ restricts to $u$, and $u_{i_1\dots
i_q}(\xi)$ restricts to $u_{i_1\dots i_q}(\I x)$; moreover, $\Id$
restricts to $\Mx$.  Hence, if the $u_{i_1\dots i_q}$ vanish at $\I
x$, then $u\in \Mx^r$.  The converse holds because the $\delta_i$
restrict to the functions $\pi_i:=\xi_i-\xi_i(\I x)$, which form a
system of regular parameters at $\I x$, and so the various products
$\pi_{i_1}\cdots\pi_{i_q}$ form a vector-space basis of $\Mx^q/
\Mx^{q+1}$.

Say that the members  $U_{\I t}$ of the linear system are of the form,
        $$U_{\I t}:t_0u^{(0)}+\cdots +t_su^{(s)}=0,$$
 where $\I t:=(t_0,\dots,t_s)\in\I P^s$ and where
$u^{(0)},\dots,u^{(s)}$ are now (regular) functions on $X$.  Expand each
$u^{(i)}$ as above, and let $u^{(i)}_{i_1\dots i_q}(\xi)$ be the
coefficients.  Then the set $M_r$ of various $r$-fold points is defined
in $X\x\I P^s$ by the vanishing of all the functions,
$$t_0u^{(0)}_{i_1\dots i_q}(\xi)+\cdots+t_su^{(s)}_{i_1\dots i_q}(\xi),$$
 where now $u(\xi)$ means the pullback of the function $u$ on $X$.  Thus
$M_r$ is closed.

Since the given point $\I x$ is a separable variable point, it lies in a
closed subset $M$ of $M_r$ such that the projection $\mu\:M\to\I P^s$ is
smooth on a dense open set.  Arguing as above in our version of
Bertini's original proof, using the closedness of $M_r$, we may assume
that $\mu$ is smooth at
        $$\I x\in U_{\I t_0}:u^{(0)}=0.$$
 We have to prove that $u^{(k)}\in\Mx^{r-1}$ for each $k$.  As before,
fix $k$ and use power series to parameterize a variable point $\I x(t)\in
M$ such that
           $$\I x(0)=\I x\hbox{ and }\mu(\I x(t))=(1,0,\dots,0,t,0,\dots,0).$$

Since $\I x(t)\in M_r$, the power series,
  $$u_{i_1\dots i_q}^{(0)}(\I x(t))+tu_{i_1\dots i_q}^{(k)}(\I x(t)),$$
vanishes identically.  Now $\Id$ is generated over the ring of
functions on $X$ by the differences $g(\eta)-g(\xi)$ as $g$ runs
through the functions on $X$.  Hence, thanks to the expansion of
$u^{(0)}(\eta)+tu^{(k)}(\eta)$, we can find functions $g_i$ and
$g_{ij}$ on $X$ such that
    $$u^{(0)}+tu^{(k)}=\sum_i g_i\big(g_{i1}-g_{i1}(\I x(t))\big)\cdots
                \big(g_{ir}-g_{ir}(\I x(t))\big).$$
  Now, $g_{ij}(\I x(t))=g_{ij}(\I x)+tg_{ij}'(t)$ for some power series
$g_{ij}'(t)$.  Plugging in and collecting the coefficient of $t$, we
conclude that $u^{(k)}\in\Mx^{r-1}$ because $g_{ij}-g_{ij}(\I x)\in\Mx$.
Thus Theorem~(4.2) is proved.

\sct The second theorem

In this section, we'll discuss Bertini's second fundamental theorem.  It
characterizes those linear systems whose members are all reducible (or
equivalently, whose general members are reducible), the so-called {\it
reducible linear systems}.  The theorem asserts, over the complex
numbers, that, if a reducible system has no fixed components, then its
members are sums of members of a pencil.  The pencil is linear if the
ambient space is the projective $n$-space $\I P^n$, the only case that
Bertini considered; see \cite{Bert82} and \cite{BertI, \p.231}.
However, in other cases, the pencil can be a {\it nonlinear\/}
1-parameter algebraic system.

Bertini's theorem was extended to an arbitrary ambient surface by
Enriques in this pioneering 1893 paper \cite{En93, \p.41}.  He
gave a sketchy argument; it is somewhat similar to the one that he gave
at the same time for Theorem~(4.1), which is described briefly in
Section~4.  In a footnote, Enriques suggests that his work be compared
to Noether's work, \cite{Noe73, \p.171} and \cite{Noe75, \p.524}, which
appeared in 1873 and 1875, about ten years before Bertini's original
paper \cite{Bert82}.  However, Noether did not, in fact, clearly state
or prove any particular case of the theorem.  Moreover, Bertini did not
refer to Noether's work in either his paper \cite{Bert82} or in his book
\cite{BertI}.  At any rate, the theorem was called the theorem of
Bertini--Enriques by van der Waerden in his 1937 paper\cite{vdW37, \S3},
and by Zariski in his 1941 paper \cite{Zar41}.  However, seventeen years
later, when Zariski \cite{Zar58} published his extension of Bertini's
theorem to positive characteristic, he did not cite Enriques' paper, nor
even mention his name.

The papers of van der Waerden and of Zariski are devoted to providing
rigorous proofs of Bertini's second theorem for an arbitrary ambient
variety of any dimension.  Bertini, Enriques (in his later treatment
\cite{EnCa32, \pp.31--33}), and van der Waerden approached the theorem
in roughly the same way, involving Bertini's first theorem.  Zariski
changed this approach at several places; he introduced some serious
commutative algebra, which in particular eliminated the need for
Bertini's first theorem.  These approaches will now be examined and
developed, but first we'll discuss the content of the second theorem
itself.

Let $X$ be an arbitrary (irreducible) variety of dimension $n$ over the
complex numbers.  A trivial way to construct a reducible linear system
on $X$ is to add a fixed component to every member of a given system.  A
more sophisticated way is to begin with a {\it pencil} without fixed
components.  This is an algebraic system of divisorial cycles such that
precisely one member passes through a given general point $\I x$ of $X$;
also, the total space $T$ is assumed to be reduced and irreducible, and
the parameter space $C$, to be complete.  If this $\I x$ lies outside
the singular locus of $X$, then $\I x$ lies inside the open set $V$ over
which the projection $T\to X$ is an isomorphism (by Zariski's ``main
theorem'').  Thus $V$ is nonempty.  Hence the total space $T$ is of
dimension $n$, and the parameter space $C$ is a curve.  So if a point
$\I y$ of $X$ lies outside the base locus of the pencil (the
intersection of all its members), then $\I y$ lies in at most finitely
many members.  Hence, if $\I y$ also lies outside the singular locus of
$X$, then $\I y$ lies in $V$.

There is a natural map $f\:V\to C$; its graph is the preimage of $V$ in
the total space $T$.  The pencil consists of the closures of the map's
fibers, aside possibly from finitely many members that have, as
components, parts of the singular locus of $X$.  Now, a linear system
(or {\it linear involution\/}) on $C$ gives rise, via pullback, to one
on $X$; its members are composed of $d$ members of the pencil if the
system on $C$ is of degree $d$.  Such a linear system on $X$ is said to
be {\it composite with a pencil\/} (or an {\it involution in a
pencil\/}).

Thus we have constructed two sorts of reducible linear systems.
Remarkably, these are the only possibilities, according to Bertini's
second fundamental theorem.  Bertini himself treated only the case where
$X=\I P^n$.  In this case, the pencil in question is necessarily a
linear system, and his theorem may be stated as follows.
        \proclaim Theorem 1 \(On reducible linear systems).
  On $\I P^n$, a reducible linear system, without fixed components, is
necessarily composite with a linear pencil.
        \endproclaim\noindent
 Bertini also restated the theorem algebraically essentially as follows.
        \proclaim Theorem 2 \(Algebraic restatement).  Let
$u^{(0)},\dots,u^{(s)}$ be forms of the same degree in variables
$x_0,\dots,x_n$, and let $t_0,\dots,t_n$ be indeterminate parameters.
If the form in the $x_i$,
        $$F:=t_0u^{(0)}+\cdots+t_su^{(s)},$$
  is reducible, then either the $u^{(i)}$ have a common factor, or they
are equal to forms $v^{(i)}$ of the same degree in two other forms $w_0$
and $w_1$ of the same degree in the $x_i$,
        $$u^{(i)}(x_0,\dots,x_n)
        = v^{(i)}(w_0(x_0,\dots,x_n),w_1(x_0,\dots,x_n)).$$
         \endproclaim

 Note that $F$ may be reducible in the $x_i$, although it is
irreducible in the $t_i$ and the $x_i$ together.  For example, $F$
might be
       $$t_0x_0^2+t_1x_1^2=\big(\sqrt{t_0}x_0+\sqrt{-t_1}x_1\big)
               \big(\sqrt{t_0}x_0-\sqrt{-t_1}x_1\big);$$
 here, $v^{(i)}(x_0,x_1)=x_i^2$ and $w^{(i)}(x_0,x_1)=x_i$.

In his book \cite{BertI, \p.231}, Bertini cited two of L\"uroth's works,
\cite{Lu93} and \cite{Lu94}, published in 1893 and 1894; in them,
Theorem~(5.2) is credited to Bertini and reproved via a lengthier, more
elementary, and more algebraic proof.  Just before, in 1891 Poincar\'e
\cite{Poin, \p.183} published a short and more algebraic proof of the
theorem for the plane, but he doesn't credit the result.  Both
Poincar\'e and L\"uroth assumed that $F$ is reducible whenever the $t_i$
are given complex values, and this hypothesis amounts precisely to the
corresponding hypothesis in Theorem~(5.1).  Bertini apparently did not
see the need to prove the equivalence of reducibility for indeterminate
values and of reducibility for all complex values; more about this gap
will be said below.

Bertini's proof runs, more or less, as follows.  By Theorem~(2.2), the
system has no variable singular points outside of the base locus.  In
particular, a general member $U$ is reduced; that is, its components
$U_1,\dots,U_d$ are distinct.  Among such $U$, fix one with $d$
minimal.  By hypothesis, $d\ge2$.  Let $\C P$ be any subsystem with the
following three properties: (1)~$\C P$ contains $U$; (2)~$\C P$ has no fixed
components; and (3)~$\C P$ is parameterized by a line $D$ in the parameter
projective space $\I P^s$ of the given linear system.  Then a general
member of $\C P$ has $d$ components because of the minimality of $d$.
Furthermore, the $i$th component varies in a 1-parameter algebraic
system $\C P_i$ , which may be constructed as follows.

Over a suitable finite algebraic extension field of the field of
rational functions on $D$, factor the form defining the total space
of $\C P$.  Say that the $i$th factor has degree $r_i$, and set
$s_i:={r_i+n\choose n}-1$.  Then the factor's coefficients are the
coordinates of a ``generic'' point of a curve $D_i$ in $\I P^{s_i}$,
which parameterizes the system $\C P_i$ of hypersurfaces of degree
$r_i$ in $\I P^n$.  In fact, $\C P_i$ is a linear pencil, because
$\deg D_i$ is equal to the number of hypersurfaces belonging to $\C
P_i$ and passing through a general point $\I x$ of $\I P^n$, and this
number is $1$ as only one member of $\C P$ passes through $\I x$.

The $d$ systems $\C P_i$ are, in fact, equal to each other.  Indeed,
given a general member $U_1'$ of $\C P_1$, say $U'_1$ is a component of
$U'\in\C P$, and choose a general point $\I x$ in $U_1'$.  For each $i$,
necessarily $\I x$ lies in some member $U'_i\in\C P_i$.  Then $U'_i$ must
also be a component of $U'$ because $\C P$ is a pencil.  Hence $U'_i$
must be equal to $U_1'$.  Thus $U_1'\in\C P_i$.  Hence $\C P_1$ and $\C
P_i$ coincide for all $i$, as asserted.  Therefore, all the components
of $U$ and of $U'$ belong to $\C P_1$.  In particular, $U_1$ and $U_2$
do; so they determine this linear pencil.  Hence, when $\C P$ is varied,
$\C P_1$ remains fixed.  Therefore, the given system is composite with
$\C P_1$.  Thus Theorem~(5.1) is proved.

Bertini gave this proof, more or less, in both \cite{Bert82} and
\cite{BertI}; in the latter, he gave more details.  The only significant
gap is the lack of justification for the assertion that a general member
of the pencil $\C P_i$ is irreducible, or as Bertini put it, that each
component of a general member varies in an algebraic system.  This gap
was filled by van der Waerden \cite{vdW37} using the theory of algebraic
cycles that he and his thesis student Wei-Liang Chow (1911--95) had just
developed (and which has become known as the theory of ``Chow
coordinates'').

In fact, van der Waerden filled the corresponding gap in Enriques'
extension of the Bertini's theorem.  Van der Waerden cited and followed
Enriques' treatment in \cite{EnCa32, \pp.31--33}, which is rather clean.
Whereas Bertini worked with an ambient projective space of arbitrary
dimension, Enriques worked only with a surface; however, van der
Waerden observed that the proof works equally well on an arbitrary
ambient variety.

Here is the statement of the general extension of Bertini's second
theorem.
        \proclaim Theorem 3 \(On reducible linear systems, extended).
On an arbitrary variety, a reducible linear system, without fixed
components, is necessarily composite with a pencil (which need not be
linear).
        \endproclaim

Enriques' proof \cite{EnCa32, \p.33}, as completed by van der Waerden
\cite{vdW37, \p.709}, runs, more or less, as follows.  Let $B$ be the
union the singular locus of the ambient variety $X$ and of the base
locus of the given system $\C U$.  Let $U$ be a general member, and
$U_1,\dots,U_d$ its components.  Then the $U_i$ are distinct; in fact,
they have no common point outside $B$, because any such point would be a
singular point of $U$, and contradict Bertini's first theorem,
Theorem~(4.1).

As in the proof of Theorem~(5.1), let $\C P$ be any subsystem with the
three properties stated there.  Then a general member of $\C P$ has $d$
components, and each varies in a certain pencil $\C P_1$, which is
independent of the component; this time, $\C P_1$ needn't be linear,
and can be defined using Chow coordinates.  Let $b$ be a base point of
$\C P$.  Since $\C P$ has no fixed components, $b$ must lie in
infinitely many members of $\C P_1$, so in all of them.  So $b$ must be
common to all the components $U_i$ of $U$.  Hence, by the conclusion of
the preceding paragraph, $b\in B$.

Consider the map $f\:(X-B)\to \I P^s$ defined by $\C U$.  Let
$\I x\in(X-B)$.  Then
        $$f^{-1}f\I x = \bigcap\{U\in\C U|\I x\in U\}.$$
  Suppose $\I x$ is general.  Then $f^{-1}f\I x$ is equidimensional.  Its
codimension must be 1; otherwise, since the $U\in \C U$ that contain
$\I x$ form a linear system, $\C U$ would contain a linear pencil having
$\I x$ as a base point, but no fixed component, contradicting the
preceding argument (applied to this pencil and a general member $U$).
It follows that, for any choice of $\C P$, necessarily $\C P_1$ is the
natural pencil $\C N$ parameterized by the normalization of $Y$ in the
function field of $X$.  Therefore, $\C U$ is composite with this
pencil, and the theorem is proved.

Enriques and van der Waerden carried out the proof above in two steps,
obtaining the following intermediate results about a linear system
without fixed components:\smallbreak
  \item{(1)} {\it If the system is reducible, then any two general
members have only base points in common.}
  \item{(2)} {\it If any two general members have only base points in
common, then the system is composite with a pencil (and conversely)}.
        \smallbreak
 On a surface, a linear system has an important invariant, its {\it
degree}; by definition, this is the number of points, other than the
base points, that two general members have in common.  In terms of this
notion, Enriques \cite{EnCa32, \p.32} restated (2) as follows: \smallbreak
 \item{(2$'$)} {\it On a surface, a linear system, without fixed
components, that has degree zero is composite with a pencil (and
conversely)}.\smallbreak\noindent
  For the plane, this result already had been proved by Bertini,
\cite{Bert01} and \cite{BertI, \p.234}.  In fact, Bertini defined the
degree of a linear system on $\I P^r$ for $r\ge3$ as well, and he
extended the characterization of those systems of degree zero, but
this extension is no longer equivalent to (2).

Zariski did not use the theory of Chow coordinates, but proved the
necessary irreducibility results directly.  In particular, he gave a
direct proof of the irreducibility of a general member of the pencil $\C
N$.  In fact, he observed that the proof works for any pencil, linear or
not, yielding the following result \cite{Zar41, \p.61}:
 \smallbreak
 \item{}{\it If a pencil, without fixed components, is not composite,
then all but finitely many members are irreducible.}\smallbreak
 Zariski also eliminated the use of Bertini's first theorem by giving
a direct algebraic proof that the pencils $\C P_1$ and $\C N$
coincide.  In the middle of this proof, he explicitly used the
hypothesis of characteristic 0.  Later, in his 1958 monograph
\cite{Zar58, \S I.6}, he treated the case of an algebraically closed
ground field of arbitrary characteristic $p$, proving this result:
 \smallbreak
 \item{}{\it If a reducible linear system is free of fixed components
and is not composite, then there is a power $p^e$ such that its
members are all of the form $p^eU$ where $U$ moves in an irreducible
linear system.}\smallbreak
 Another proof had already been published in 1951 by Matsusaka
\cite{Mat51}.  However, the author was a student  of Zariski's in the
early 1960s, and can distinctly recall several conversations with him in
which he shook his head and said: ``Poor Matsusaka, he didn't know that
I already had a proof.''

\sctvspace
\centerline{REFERENCES}
\references

Bert77
 E. Bertini,
 {\it Ricerche sulle trasformazioni univoche involutorie nel piano,}
 \adm \(2)8 1877 244-86

Bert82
 E. Bertini,
 {\it Sui sistemi lineari,}
 \rril \(II)15 1882 24--28

Bert01
 E. Bertini,
 {\it Sui sistemi lineari di grado zero,}
 \anlr \(5)10 1901 73--76

Bert04
 E. Bertini,
 {\it Life and works of L. Cremona},
 Proc. Lond. Math. Soc., (2){\bf 1} (1904), V--XVIII.
 Edited and expanded as ``Della vita e delle opere di L. Cremona'' in
 \vol.3 of  {\it Opere matematiche di\/} {\smc Lugi Cremona},
 Milano, 1917, V--XXII.

BertI
 E. Bertini,
 ``Introduzione alla geometria proiettiva degli iperspazi,''
 Enrico Spoerri, Pisa, 1907.

BertC
 E. Bertini,
 ``Complementi di geometri proiettiva,''
 Bologna, 1928.

Berz15
 L. Berzolari,
 {\it Th\'eorie g\'en\'erale des courbes planes alg\'ebriques,}
 in ``Encyclop\'edie des sciences math\'ematiques, III.19,''
 Gauthier-Villars, Paris (1915), 257--304.

Berz33a
 L. Berzolari,
 {\it Eugenio Bertini,}
 \rril\(2)66 1933 67--68, 609--35

Berz33b
 L. Berzolari,
 {\it Eugenio Bertini,}
 \bumi  12 1933 148--53

Bot77
 U. Bottazzini,
 {\it Riemanns Einfluss auf E Betti und F Casorati,}
 \ahes 18\(1) 1977/78 27--37

Bot81
 U. Bottazzini,
 {\it Mathematics in a unified Italy,}
 in ``Social history of nineteenth century mathematics,''
 Mehrtens, Bos, Schneider (eds.),
 Birkha\"user (1981), 165--78.

BoCoGa96
 U. Bottazzini, A. Conte and P. Gario,
 ``Riposte Armonie,''
 Bollati Boringhieri, To\-ri\-no, 1996.

BrCi95
 A. Brigaglia and C. Ciliberto,
 ``Italian Algebraic Geometry between the Two World Wars,''
 Queen's Papers in Pure and Applied Math, \vol. 100, 1995.

BrCiSe92
 A. Brigaglia{,} C. Ciliberto  and E. Sernesi{,} ed.,
 ``Algebra e Geometria (1860--1940): Il contributo italiano,''
 Conference Proc., Cortona, Italy, 4--8 May 1992,
 Supplemento Rend. Circ. Mat. Palermo, (II){\bf 36} 1994.

Bri92
 F. Brioschi,
 {\it Enrico Betti,}
 \adm \(II)20 1892 256 ($=$  {\it Opere matematiche,} \vol.3, Hoepli,
Milano, 1904)

Cast33
 G. Castelnuovo,
 {\it Commemorazione del socio Eugenio Bertini,}
 \anlr \(6)17 1933 745--48

Con34
 A. Conti,
  {\it Onoranze ad Eugenio Bertini},
 \bdmf \(2)13 1934 46--52

Cool
 J. L. Coolidge,
 ``A treatise on algebraic plane curves,''
 Oxford, 1931 ($=$ Dover reprint, 1959).

Dieu
 J. Dieudonn\'e,
 ``Cours de G\'eom\'etrie alg\'ebrique, \vol.1 (Aper\c cu historique sur
le d\'eveloppement de la g\'eom\'etrie alg\'ebrique),''
Press. Univ. France, 1974.

En93
 F. Enriques,
 {\it Richerche di geometria sulle superficie algebriche,}
 \mast \(2)44 1983 171--232 ($=$\ms I 31--106)

En96
 F. Enriques,
 {\it Introduzione alla geometria sopra le superficie algebriche,}
 \xlm \(3)10 1896 1--81 ($=$\ms I 211--312)

En38
 F. Enriques,
 ``Le matematiche nella storia e nella cultura,''
 Zanichelli, Bologna, 1938 ($=$ reprint, 1982).

EnCa32
 F. Enriques and L. Campedelli,
 ``Lezioni sulla teoria delle superficie algebriche,''
 Lithographed, CEDAM, Padova, 1932.

EnCh18
 F. Enriques and O. Chisini,
 ``Lezioni sulla teoria geometrica delle equazioni e delle funzioni
algebriche, \vol.1, 2''
 Zanichelli, Bologna, 1915, 1918.

Fub33
 G. Fubini,
 {\it Eugenio Bertini,}
 \aast 68 1933 447--53

Greit
 S. Greitzer,
 {\it Antonio Luigi Gaudenzio Giuseppe Cremona,}
 in ``Dictionary of Scientific Biography,''
 Gillispie (ed.), \vol. III, Scribner, 1971.

Lu93
 J. L\"uroth,
 {\it Bewis eines Satzes von Bertini \"uber lineare Systeme ganzer
Funtionen,}
 \ma 42 1893 457--70

Lu94
 J. L\"uroth,
 {\it Bewis eines Satzes von Bertini \"uber lineare Systeme ganzer
Funtionen. II,}
 \ma 44 1894 539--52

Mat51
 T. Matsusaka,
 {\it The theorem of Bertini on linear systems in modular fields,}
 \mcsuk 26 1951 51--62

Mum72
 D. Mumford,
 {\it Introduction,}
  in ``Oscar Zariski: Collected Papers I,''
 Hironaka and Mumford (ed.),
 MIT Press (1972),
 3--5.

Noe73
 M. Noether,
 {\it Ueber Fl\"achen, welche Schaaren rationaler Curven besitzen,}
 \ma 3 1873 161--227

Noe75
 M. Noether,
 {\it Zur Theorie des eindeutigen Entsprechens algebraischer Gebilde,}
 \ma 8 1875 495--533

Par91
 C. Parikh,
 ``The unreal life of Oscar Zariski,''
 Academic Press, 1991.

PiSi06
 E. Picard and G. Simart,
 ``Th\'eorie des Fonctions Alg\'ebriques des deux variables
ind\'ependantes, Tome II,''
 Gauthier-Villars, Paris, 1906.

Poin
 H. Poincar\'e,
 {\it Sur l'int\'egration alg\'ebrique des \'equations diff\'erential du
premier order et du premier degr\'e,}
 \rcmp 5 1891 161--91

Ros71
 J. Rosanes,
 {\it Ueber diejenigen rationalen Substitutionen, welche eine rationale
Umkehrung zulassen,}
 \jram 73 1871 97--100

Scor34
 G. Scorza,
 {\it Eugenio Bertini,}
 \emcmc\(2)7 1934 101-17

Sev06
 F. Severi,
 {\it Su alcuni propriet\`a dei moduli di forme algebriche,}
 \aast \(2)41 1906 167--85

Sev08
 F. Severi,
 ``Lezioni di Geometria Algebrica,''
 Angelo Draghi, Padova, 1908.

vdW35
 B. L. van der Waerden,
 {\it Zur algebraischen Geometrie. V. Ein Kriterium f\"ur die
Einfachheit von Schnittpunkten,}
 \ma 110 1935 128--33

vdW37
 B. L. van der Waerden,
 {\it Zur algebraischen Geometrie. X. Uber linear Scharen von
reduziblen Mannigfaltigkeiten,}
 \ma 113 1937 705--12

Ver03
 G. Veronese,
 {\it Commemorazione del socio Luigi Cremona,}
 \rral \(5)12 1903 664--78

Volt
 V. Volterra
 {\it Betti, Brioschi, Casorati, trois analystes italiens et trois
mani\`eres d'envisager les questions d'analyse,} in ``Compte Rendu
deuxi\`eme Congr. des Math., Gauthier-Villars, Paris, 1902.

Zar35
 O. Zariski,
 ``Algebraic Surfaces,''
 Ergebnisse Math., Band II, Heft 5, Springer-Verlag 1935.
 Second Suppl. Ed., Ergebnisse Math., Band 61,
 Springer-Verlag 1971.

Zar41
 O. Zariski,
 {\it Pencils on an algebraic variety and a new proof of a theorem of
Bertini,}
 \tams50 1941 48--70 ($=$ \cp I 154--76)

Zar44
 O. Zariski,
 {\it The theorem of Bertini on the variable singular points of a linear
system of varities,}
 \tams 56 1944 130--40 ($=$ \cp I 242--52)

Zar44r
 O. Zariski,
 {\it Reduction of the singularities of algebraic three dimensional
varieties,}
 \am 45 1944 472--542 ($=$ \cp I 441--511)

Zar47
 O. Zariski,
 {\it The concept of a simple point of an abstract algebraic variety,}
 \tams 62 1947 1-52 ($=$ \cp I 252--304)

Zar58
 O. Zariski,
 ``Introduction to the problem of minimal models in the theory of
algebraic surfaces,"
 Publ. Math. Soc. Japan 4, 1958
 ($=$ \cp II 277--369).

\endreferences

\line{\hskip9.5truecm\hrulefill}\medskip
\rightline{\it Department of Mathematics}
\rightline{\it University of Copenhagen}
\rightline{\it Universitetsparken 5}
\rightline{\it DK 2100 Copenhagen \O}
\rightline{\it Denmark}
\medskip
\rightline{\it Department of Mathematics}
\rightline{\it Room 2-278 MIT}
\rightline{\it 77 Massachusetts Avenue}
\rightline{\it Cambridge, MA 02139}
\rightline{\it USA}

\bye